\pdfoutput=1
\RequirePackage{ifpdf}
\ifpdf 
\documentclass[pdftex]{sigma}
\else
\documentclass{sigma}
\fi

\usepackage{appendix}
\usepackage{booktabs}
\usepackage{empheq}
\usepackage{enumitem}
\usepackage{fancyvrb}
\usepackage{multirow}
\usepackage{pdflscape}

\newtheorem*{theoremA}{Theorem A}

\newlength{\caseIndent}
\DeclareMathOperator{\ad} {ad}
\DeclareMathOperator{\Der} {Der}
\DeclareMathOperator{\diag}{diag}
\DeclareMathOperator{\G} {G}
\DeclareMathOperator{\id} {id}
\DeclareMathOperator{\rank}{rank}
\DeclareMathOperator{\sign}{sign}
\DeclareMathOperator{\SO} {SO}
\DeclareMathOperator{\Sp} {Sp}
\DeclareMathOperator{\spanOperator}{span}
\DeclareMathOperator{\tr} {tr}

\renewcommand{\Im}{\operatorname{Im}}
\renewcommand{\Re}{\operatorname{Re}}

\newcommand{\Span}[1]{\spanOperator\{{#1}\}}
\newcommand{\wt} [1]{\smash{\widetilde{#1}}}

\AtEndDocument{\vfill\eject\batchmode} 
\usepackage{amssymb}
\usepackage[T1]{fontenc}
\DeclareSymbolFont{script}{U}{eus}{m}{n}
\DeclareMathSymbol{\Wedge}{0}{script}{"5E}

\newcommand{\bbC}{\mathbb{C}}
\newcommand{\bbF}{\mathbb{F}}
\newcommand{\bbI}{\mathbb{I}}
\newcommand{\bbJ}{\mathbb{J}}
\newcommand{\bbO}{\mathbb{O}}
\newcommand{\bbP}{\mathbb{P}}
\newcommand{\bbR}{\mathbb{R}}

\newcommand{\mbA}{\mathbf{A}}
\newcommand{\mbB}{\mathbf{B}}
\newcommand{\mbC}{\mathbf{C}}
\newcommand{\mbD}{\mathbf{D}}
\newcommand{\mbE}{\mathbf{E}}
\newcommand{\mbF}{\mathbf{F}}
\newcommand{\mbH}{\mathbf{H}}
\newcommand{\mbK}{\mathbf{K}}
\newcommand{\mbL}{\mathbf{L}}
\newcommand{\mbO}{\mathbf{O}}
\newcommand{\mbr}{\mathbf{r}}
\newcommand{\mbS}{\mathbf{S}}
\newcommand{\mbs}{\mathbf{s}}
\newcommand{\mbT}{\mathbf{T}}
\newcommand{\mbU}{\mathbf{U}}
\newcommand{\mbV}{\mathbf{V}}
\newcommand{\mbv}{\mathbf{v}}
\newcommand{\mbW}{\mathbf{W}}
\newcommand{\mbX}{\mathbf{X}}
\newcommand{\mbY}{\mathbf{Y}}
\newcommand{\mbZ}{\mathbf{Z}}
\newcommand{\mbzero}{\mathbf{0}}

\newcommand{\mcL}{\mathcal{L}}

\newcommand{\mfd} {\mathfrak{d}}
\newcommand{\mfe} {\mathfrak{e}}
\newcommand{\mff} {\mathfrak{f}}
\newcommand{\mfg} {\mathfrak{g}}
\newcommand{\mfh} {\mathfrak{h}}
\newcommand{\mfk} {\mathfrak{k}}
\newcommand{\mfm} {\mathfrak{m}}
\newcommand{\mfN} {\mathfrak{N}}
\newcommand{\mfn} {\mathfrak{n}}
\newcommand{\mfq} {\mathfrak{q}}
\newcommand{\mfr} {\mathfrak{r}}
\newcommand{\mft} {\mathfrak{t}}
\newcommand{\mfaut}{\mathfrak{aut}}
\newcommand{\mfgl} {\mathfrak{gl}}
\newcommand{\mfnil}{\mathfrak{nil}}
\newcommand{\mfrad}{\mathfrak{rad}}
\newcommand{\mfsl} {\mathfrak{sl}}
\newcommand{\mfso} {\mathfrak{so}}
\newcommand{\mfsp} {\mathfrak{sp}}

\numberwithin{equation}{section}

\newtheorem{Theorem}{Theorem}[section]

\newtheorem{Proposition}[Theorem]{Proposition}

 { \theoremstyle{definition}
\newtheorem{Remark}[Theorem]{Remark}
\newtheorem{Example}[Theorem]{Example}}

\begin{document}

\newcommand{\arXivNumber}{1807.02734}

\renewcommand{\PaperNumber}{008}

\FirstPageHeading

\ShortArticleName{Homogeneous Real $(2,3,5)$ Distributions with Isotropy}

\ArticleName{Homogeneous Real $\boldsymbol{(2,3,5)}$ Distributions with Isotropy}

\Author{Travis WILLSE}

\AuthorNameForHeading{T.~Willse}

\Address{Fakult\"{a}t f\"{u}r Mathematik, Universit\"{a}t Wien, Oskar-Morgenstern-Platz 1, 1090 Wien, Austria}
\Email{\href{mailto:travis.willse@univie.ac.at}{travis.willse@univie.ac.at}}

\ArticleDates{Received August 15, 2018, in final form January 26, 2019; Published online February 04, 2019}

\Abstract{We classify multiply transitive homogeneous real $(2,3,5)$ distributions up to local diffeomorphism equivalence.}


\Classification{53A30; 53C15; 53C30}

\section{Introduction}

The study of $(2, 3, 5)$ distributions, that is, tangent $2$-plane distributions $(M, \mbD)$ on $5$-manifolds satisfying the genericity condition $[\mbD, [\mbD, \mbD]] = TM$, dates to Cartan's celebrated ``Five Variables Paper'' \cite{CartanFiveVariables}. That article (1) resolved the equivalence problem for this geometry, (2) in doing so revealed a surprising connection with the exceptional complex Lie algebra of type $\G_2$, and (3) (nearly) locally classified complex $(2, 3, 5)$ distributions $\mbD$ whose infinitesimal symmetry algebra $\mfaut(\mbD)$ has dimension at least $6$. Besides its historical significance and its connection with $\G_2$, which mediates its relationship with other geometries \cite{CapSagerschnig, GPW, GrahamWillse, HammerlSagerschnig, HammerlSagerschnigTwistor, LeistnerNurowski, LNS}, \cite[Section~5]{Nurowski}, this geometry is significant because of its appearance in simple, nonholonomic kinematic systems \cite{BorMontgomery, BryantHsu}. It has enjoyed heightened attention in the last decade or so \cite{AnNurowski, AnNurowskiPlebanskiMetric, DaveHaller, DoubrovGovorovCounterexample, DHS, DoubrovKruglikov, IKTY, IKY, Randall,SagerschnigWillseOperator, SagerschnigWillseGeometry, Strazzullo, WillseHeisenberg, WillseCounterexample}, owing in part to its realization in the class of parabolic geometries \cite[Section~4.3.2]{CapSlovak}, \cite{Sagerschnig, SagerschnigThesis}, a broad family of Cartan geometries for which many powerful results are available.

In the current article we classify locally all homogeneous \textit{real} $(2, 3, 5)$ distributions with multiply transitive symmetry algebra, so again those for which $\dim \mfaut(\mbD) \geq 6$. Our motivation for carrying out this classification is twofold: (1) It gives a canonical list of examples with favorable symmetry properties; this list is exploited in work in progress about real $(2, 3, 5)$ distributions. (2) It is independently interesting, in part because of the appearance of several distinguished rolling distributions; see Section~\ref{section:rolling}.

Our method is standard: Any homogenous $(2, 3, 5)$ distribution $(M, \mbD)$ can be encoded in an \textit{algebraic model} $(\mfh, \mfk; \mfd)$ in the sense that the original distribution can be recovered (up to local diffeomorphism equivalence) from and hence specified by the model data. Here, $\mfh := \mfaut(\mbD)$ is the (infinitesimal) symmetry algebra of $\mbD$, $\mfk$ is the isotropy subalgebra of a point $u \in M$ (so, in our multiply transitive setting, $\dim \mfk = \dim \geq 1$), and $\mfd \subset \mfh$ is the subspace corresponding to $\mbD_u$. Given any real algebraic model, its complexification $(\mfh \otimes \bbC, \mfk \otimes \bbC; \mfd \otimes \bbC)$ is a complex algebraic model, and conversely the given real algebraic model can be recovered from an appropriate anti-involution $\phi\colon \mfh \otimes \bbC \to \mfh \otimes \bbC$ admissible in the sense that it preserves the filtration $\mfk \otimes \bbC \subset \mfd \otimes \bbC \subset \mfh \otimes \bbC$.

We thus briefly summarize in Section~\ref{section:classification-complex} and record in Tables~\ref{table:complex-14}--\ref{table:complex-6} in Appendix~\ref{appendixB} the classification of complex $(2, 3, 5)$ distributions. For each distribution in the classification we give an explicit algebraic model in terms of abstract Lie algebra data. Most of these distributions were identified by Cartan~\cite{CartanFiveVariables}, but Doubrov and Govorov found much later that Cartan's list omitted the model we here call~$\mathbf{N.6}$~\cite{DoubrovGovorovCounterexample}.

In Section~\ref{section:classification-real} we classify for each complex algebraic model $(\mfh, \mfk; \mfd)$ with $\dim \mfh = 6$ (these are recorded in Table~\ref{table:complex-6}) the admissible anti-involutions of $\mfh$ up to a notion of equivalence that corresponds to diffeomorphism equivalence of the homogeneous distributions and record the corresponding fixed real Lie algebra data. Together with the real algebraic model $\mbO^{\bbR}$ with maximal symmetry (unique up to equivalence) and the pre-existing classification of the models with $\dim \mfh = 7$ \cite[Theorem~2]{DoubrovZelenko}, here denoted $\mathbf{N.7}_{\Lambda}^{\bullet}$, this yields the main result of this article:

\begin{theoremA}\label{theorem:classification-real} Any multiply transitive homogeneous real $(2, 3, 5)$ distribution is locally equivalent to exactly one distribution in Tables~{\rm \ref{table:real-O}--\ref{table:real-D.6-ast}}, modulo the given equivalences of parameters.
\end{theoremA}

In Section~\ref{section:identification} we give algorithms for identifying, in both the complex and real cases, a multiply transitive homogeneous $(2, 3, 5)$ distribution given in terms of an abstract algebraic model among the distributions in the classification; this amounts to constructing sufficiently many invariants to distinguish all of the models. Finally, in Section~\ref{section:rolling} we identify many of the distributions in the real classification as rolling distributions, that is, $2$-plane distributions on $5$-manifolds defined on the configuration space of two surfaces rolling along one another by the kinematic no-slip and no-twist conditions.

\section[The geometry of $(2,3,5)$ distributions]{The geometry of $\boldsymbol{(2,3,5)}$ distributions}

A \textit{$(2, 3, 5)$ distribution} is a tangent $2$-plane distribution $\mbD$ on a $5$-manifold $M$ satisfying the genericity condition
\begin{gather*}
	[\mbD, [\mbD, \mbD]] = TM .
\end{gather*}
Here, for distributions $\mbE, \mbE'$ on $M$, $[\mbE, \mbE'] := \coprod_{u \in M} [\mbE, \mbE']_u$, where $[\mbE, \mbE']_u \subset T_u M$ is the vector subspace $\{[X, Y]_u \colon X \in \Gamma(\mbE), Y \in \Gamma(\mbE')\}$. Implicit in the notation $[\mbD, [\mbD, \mbD]]$ is the condition that $[\mbD, \mbD]$ has constant rank; for a $(2, 3, 5)$ distribution, $\rank{}[\mbD, \mbD] = 3$. We will work in both the smooth and complex categories. In both cases we will always assume that~$M$ is connected.

\subsection{Monge (quasi-)normal form}

Any ordinary differential equations of the form
\begin{gather}\label{equation:ode}
	z'(x) = F(x, y, y', y'', z)
\end{gather}
defines a $2$-plane distribution on the (partial) jet space $J^{2, 0}\big(\bbF, \bbF^2\big) \cong \bbF^5$ ($\bbF = \bbR$ or $\bbF = \bbC$). We can prolong any solution $(x, y(x), z(x))$ to a curve $(x, y(x), y'(x), y''(x), z'(x))$ in $\bbF^5$, and by construction any such curve is tangent to the $2$-plane distribution $\mbD_F \subset T \bbF^5$ defined in the respective jet coordinates $(x, y, p, q, z)$ as the common kernel of the canonical jet $1$-forms ${\rm d}y - p \,{\rm d}x$ and ${\rm d}p - q \,{\rm d}x$ and the $1$-form ${\rm d}z - F(x, y, p, q, z) \,{\rm d}x$ defined by~\eqref{equation:ode}. Conversely, the projection to $xyz$-space of any integral curve of this distribution to which the pullback of $dx$ vanishes nowhere defines a solution of this o.d.e. The distribution $\mbD_F$ is spanned by
\begin{gather*}
	\partial_q
		\qquad \textrm{and} \qquad
	D_x := \partial_x + p \partial_y + q \partial_p + F(x, y, p, q, z) \partial_z
\end{gather*}
--~the latter is the total derivative~-- and computing directly shows that $\mbD_F$ is a $(2, 3, 5)$ distribution iff $\partial_q^2 F$ vanishes nowhere.

Goursat showed that, in fact, every $(2, 3, 5)$ distribution arises locally this way and hence can be specified locally by some function $F$ of five variables. Such an o.d.e.\ (or, by slight abuse of terminology, the function $F$ itself) is called a \textit{Monge normal form} of the distribution.
\begin{Proposition}[{\cite[Section~76]{Goursat}}]
Let $(M, \mbD)$ be a real (complex) $(2, 3, 5)$ distribution and fix a~point $u \in M$. There is a neighborhood $U \subseteq M$ of $u$, a diffeomorphism $($biholomorphism$)$ $\Psi\colon U \to V \subset J^{2, 0}(\bbF, \bbF^2)$, and a smooth (complex-analytic) function $F\colon V \to \bbF$ for which $T\Psi \cdot \mbD\vert_U = \mbD_F$.
\end{Proposition}

\section{Homogeneous distributions}

\subsection{Infinitesimal symmetries} An \textit{infinitesimal symmetry} of a $(2, 3, 5)$ distribution $(M, \mbD)$ is a vector field $\xi \in \Gamma(\mbD)$ 	whose (local) flow preserves $\mbD$, or equivalently for which $\mcL_{\xi} \eta \in \Gamma(\mbD)$ for all $\eta \in \Gamma(\mbD)$. We denote the Lie algebra of infinitesimal symmetries, called the \textit{$($infinitesimal$)$ symmetry algebra}, by $\mfaut(\mbD)$, and we say that $(M, \mbD)$ is \textit{infinitesimally homogeneous} if $\mfaut(\mbD)$ acts infinitesimally transitively, that is, if $\{\xi_u \colon \xi \in \mfaut(\mbD)\} = T_u M$ for all $u \in M$. This article concerns (infinitesimally) \textit{multiply transitive} homogeneous distributions, that is, infinitesimally homogeneous distributions for which the isotropy subalgebra $\mfk_u < \mfaut(\mbD)$ of infinitesimal symmetries vanishing at any $u \in M$ is nontrivial, or equivalently, for which $\dim \mfaut(\mbD) \geq 6$.

\subsection{Algebraic models for homogeneous distributions}

Fix a homogeneous $(2, 3, 5)$ (real or complex) distribution $(M, \mbD)$ with transitive symmetry algebra $\mfh := \mfaut(\mbD)$, fix a point $u \in M$, and denote by $\mfk < \mfh$ the subalgebra of vector fields in $\mfh$ vanishing at $u$ and by $\mfd \subset \mfh$ the subspace $\mfd := \{\xi \in \mfh \colon \xi_u \in \mbD_u\}$. Then, $\mfk \subseteq \mfd$, $[\mfk, \mfd] \subseteq \mfd$ (we call this property \textit{$\mfk$-invariance}), and $\dim(\mfd / \mfk) = 2$. The fact that $\mbD$ is a $(2, 3, 5)$ distribution implies the genericity condition $\mfd + [\mfd, \mfd] + [\mfd, [\mfd, \mfd]] = \mfh$. We call the triple $(\mfh, \mfk; \mfd)$ a (real or complex) \textit{algebraic model} (of $(M, \mbD)$).

Given an algebraic model, we can reconstruct $\mbD$ up to local equivalence: For any groups $H$, $K$ with $K < H$ and respectively realizing $\mfh$, $\mfk$ (for the groups that occur in the classification, we can choose $K$ to be a closed subgroup of $H$), invoke the canonical identification $T_{\id \cdot K} (H / K) \cong \mfh / \mfk$ to take $\mbD \subset T(H / K)$ to be the distribution with fibers $\mbD_{h \cdot K} = T_{\id \cdot K} L_h \cdot (\mfd / \mfk)$, where $L_h \colon H / K \to H / K$ is the map $L_h \colon g \cdot K \mapsto (hg) \cdot K$; by $\mfk$-invariance this definition is independent of the coset representative $h$, and by genericity $\mbD$ is an $H$-invariant $(2, 3, 5)$ distribution. Via the above identification, $[\mbD, \mbD]_{\id \cdot K} = (\mfd + [\mfd, \mfd]) / \mfk$.\footnote{Conversely, a triple $(\mfh, \mfk; \mfd)$, where $\mfh$ is a Lie algebra, $\mfk < \mfh$ is a Lie subalgebra, and $\mfd \subset \mfh$ is a subspace for which $\mfh \supset \mfk$ and $\dim (\mfd / \mfk) = 2$, together satisfying the $\mfk$-invariance and genericity conditions, determines up to local equivalence a homogeneous distribution via this construction. The symmetry algebra of this distribution may be strictly larger that $\mfh$, however; for example, for $\mathbf{D.6}_{\lambda}$ we have $\dim \mfh = 6$, but for the excluded value $\lambda = 9$, the resulting distribution $\Delta$ is flat, so $\dim \mfaut(\Delta) = \dim \mfg_2 = 14$.}

We declare two algebraic models $(\mfh, \mfk; \mfd), (\mfh', \mfk'; \mfd')$ to be \textit{equivalent} iff there is a Lie algebra isomorphism $\alpha\colon \mfh \to \mfh'$ satisfying $\alpha(\mfk) = \mfk'$ and $\alpha(\mfd) = \mfd'$. Unwinding definitions shows that equivalent algebraic models determine locally equivalent distributions.

\section{Multiply transitive homogeneous complex distributions}\label{section:classification-complex}

Cartan showed that for all $(2, 3, 5)$ distributions $\mbD$, $\dim \mfaut(\mbD) \leq 14$, and that equality holds iff it is locally equivalent to the so-called flat distribution $\Delta$ \cite{CartanFiveVariables}; his argument applies to both the real and complex settings. We call the corresponding (complex) algebraic model $\mbO$ (see Section~\ref{subsection:O}). In this case, $\mfaut(\mbD)$ is isomorphic to the simple complex Lie algebra of type $\G_2$~-- we denote it by $\mfg_2(\bbC)$~-- and we say that $\mbD$ is \textit{$($locally$)$ flat}. Cartan furthermore claimed to classify up to local equivalence (and implicitly in the complex setting) all distributions $\mbD$ with $\dim \mfaut(\mbD) \geq 6$.\footnote{Cartan's classification is restricted to distributions for which the Petrov type of the distribution is the same at all points; this condition holds automatically for locally homogeneous distributions. See Remark~\ref{remark:labeling}.\label{footnote:constant-type}} Doubrov and Govorov found (much) later, however, that Cartan's classification misses a single distribution up to local equivalence (we denote it~$\mathbf{N.6}$) \cite{DoubrovGovorovCounterexample, WillseCounterexample}.

In this section, we briefly summarize the classification of multiply transitive homogeneous complex $(2, 3, 5)$ distributions $\mbD$ up to local equivalence.

\begin{Theorem}[\cite{CartanFiveVariables, DoubrovGovorovCounterexample}]\label{theorem:complex-classification}
Any multiply transitive homogeneous complex $(2, 3, 5)$ distribution is locally equivalent to exactly one distribution in Tables~{\rm \ref{table:complex-14}--\ref{table:complex-6}}, modulo the given equivalences of parameters.
\end{Theorem}

In these tables we record for each such distribution (1) an explicit algebraic model $(\mfh, \mfk; \mfd)$ in terms of abstract Lie algebra data, (2) a Monge normal form $F$ realizing the distribution, (3) an explicit isomorphism $\mfh \cong \mfaut(\mbD_F)$, and (4) a basis of $\mfh$ adapted to the algebraic model, which we exploit in the real classification.

In Section~\ref{section:classification-real} we use this list of complex algebraic models to classify the real algebraic models.

We use the convention that the undecorated Fraktur names $\mfg_2$, $\mfgl_m$, $\mfsl_m$, $\mfso_m$, $\mfsp_m$ refer to real Lie algebras, and we denote their complexifications by $\mfgl_m(\bbC)$ and analogously. By mild abuse of notation, for any element $\mbv$ of a real Lie algebra $\mfg$ we also denote by $\mbv$ the element $\mbv \otimes 1 \in \mfg \otimes \bbC$.

\begin{Remark}\label{remark:labeling}Our convention for labeling the nonflat distributions $\mbD$ in the classification refers both to the dimension of $\mfaut(\mbD)$ and to a particular discrete invariant. The fundamental curvature quantity of a $(2, 3, 5)$ distribution $(M, \mbD)$ is a section $A \in \Gamma(S^4 \mbD^*)$ \cite[Section~33]{CartanFiveVariables}, and its nonvanishing is a complete local obstruction to local equivalence to the model $\mbO$ (Section~\ref{subsection:O}) \cite[Sections~36--39]{CartanFiveVariables}, hence the epithet \textit{flat}. At each point $u \in M$ the \textit{Petrov type} (or \textit{root type}) of $A_u$ is the multiplicities of the roots of $A_u$; if $\mbD$ is real, we instead use the multiplicities of the roots of $A_u \otimes \bbC$. If $\mbD$ is infinitesimally homogeneous, the Petrov type is the same at all points, and among multiply transitive homogeneous distributions only Petrov types~D (two double roots), N (a quadruple root), and O ($A = 0$) occur.
\end{Remark}

\subsection[The flat model $\mbO$]{The flat model $\boldsymbol{\mbO}$}\label{subsection:O}
Denote by $\mfg_2$ the split real form of $\mfg_2(\bbC)$, take $\mfq < \mfg_2$ to be the (parabolic) subalgebra of elements fixing an isotropic line in the standard representation of $\mfg_2$ (cf.\ \cite[Section~4]{BorMontgomery}), and denote by $\mfq_+ < \mfq$ the orthogonal of $\mfq$ with respect to the Killing form on $\mfg_2$. Then, define the subspace $\mfg_2^{-1} := \{ \xi \in \mfg_2 \colon [\mfq_+ , \xi] \subseteq \mfq \}$. The Killing form bracket identity and the Jacobi identity together give $\big[\mfq(\bbC), \mfg_2(\bbC)^{-1}\big] \subseteq \mfg_2(\bbC)^{-1}$ (in fact, equality holds), and inspecting the root diagram of $\mfg_2$~-- or just using the explicit realization $\mfg_2 < \mfgl_7$ (see Appendix~\ref{appendix:tables})~-- gives that $\dim \big(\mfg_2^{-1} / \mfq\big) = 2$ and $\big[\mfg_2^{-1}, \big[\mfg_2^{-1}, \mfg_2^{-1}\big]\big] = \mfg_2$, so
\begin{gather*}
	\boxed{\big(\mfg_2(\bbC), \mfq(\bbC); \mfg_2^{-1}(\bbC)\big)}
\end{gather*}
is an algebraic model, $\mbO$, of the complex flat distribution.\footnote{This is a special case of a much more general construction \cite[Section~3]{CapSlovak} related to the realization of $(2, 3, 5)$ distributions as so-called parabolic geometries \cite[Section~4.3.2]{CapSlovak}, \cite{Sagerschnig}.} See Table~\ref{table:complex-14} for details.

The flat distribution can also be realized in Monge normal form by $F(x, y, p, q, z) = q^2$, or equivalently, by the so-called Hilbert--Cartan equation $z'(x) = y''(x)^2$ \cite{CartanCurves, Hilbert}.

\subsection[The submaximal models $\mathbf{N.7}_\Lambda$]{The submaximal models $\boldsymbol{\mathbf{N.7}_\Lambda}$}\label{subsection:N.7}
If $\dim \mfaut(\mbD) < 14$, then in fact $\dim \mfaut(\mbD) \leq 7$ \cite{CartanFiveVariables, KruglikovThe};\footnote{As in footnote~\ref{footnote:constant-type}, the bound in \cite{CartanFiveVariables} is established for distributions with constant Petrov type; this assumption was eliminated in \cite{KruglikovThe}.} distributions for which equality holds are sometimes called \textit{submaximal}. All submaximal complex distributions are homogeneous, and they fit together in a $1$-parameter family, for which several convenient Monge normal forms have been found, including \cite[Section~45]{CartanFiveVariables}\footnote{Here the factor $\frac{10}{3}$ corrects a numerical error of Cartan~\cite{Strazzullo}.}
\begin{gather*}
	F_I(x, y, p, q, z) := q^2 + \tfrac{10}{3} I p^2 + \big(I^2 + 1\big) y^2 , \qquad I \in \bbC .
\end{gather*}
The quantity $I^2 \in \bbC$ is a complete invariant.

It is convenient for our purposes to use a generalization of this form appearing in \cite[Section~4]{CartanFiveVariables} and studied by Doubrov and Zelenko in the context of control theory: Define \cite{DoubrovZelenko}
\begin{gather*}
	F_{r, s}(x, y, p, q, z) = q^2 + r p^2 + s y^2, \qquad r, s \in \bbC ,
\end{gather*}
and denote the distribution it determines by $\mbD_{r, s}$. Then, $\mbD_{r, s}$ is locally equivalent to the flat model iff the roots of the polynomial $t^4 - r t^2 + s$ form an arithmetic sequence, that is, if $9 r^2 = 100 s$; otherwise it is submaximal. In Section~\ref{subsubsection:identification-submaximal-complex} we recover the fact that the distributions $\mbD_{r, s}$ and $\mbD_{r', s'}$ are locally equivalent iff there is a constant $c \in \bbC - \{ 0 \}$ such that $r' = c r$, $s' = c^2 s$. For convenience we use the invariant
\begin{gather*}
 \Lambda = \frac{64 s}{100 s - 9 r^2} = \frac{16}{25} \big(I^2 + 1\big) ,
\end{gather*}
and by mild abuse of notation we denote the submaximal distribution with this invariant value by $\mbD_\Lambda$. Since $I^2$ is a complete invariant, so is $\Lambda$.

We realize the distributions $\mbD_{\Lambda}$ as algebraic models as follows. Let $\mfn = \mfn_{-2} \oplus \mfn_{-1}$ be the $5$-dimensional real Heisenberg algebra endowed with its standard contact grading, and fix a standard adapted basis $(\mbU, \mbS_1, \mbS_2, \mbT_1, \mbT_2)$ of $\mfn$, so that $\mfn_{-2} = \Span{\mbU}$, $\mfn_{-1} = \Span{\mbS_1, \mbS_2, \mbT_1, \mbT_2}$, $[\mbS_1, \mbT_1] = [\mbS_2, \mbT_2] = \mbU$, and all brackets of basis elements not determined by these identities are zero. Let $\mbE \in \Der(\mfn)$ denote the grading derivation, so that $\mbE\vert_{\mfn_\ell} = \ell \id_{\mfn_\ell}$ for $\ell = -2, -1$.

For each $\Lambda \in \bbC$, pick $(r, s)$ satisfying $\Lambda = 64 s / \big(100 s - 9 r^2\big)$ (one can always choose $r = 2 (25 \Lambda - 16), s = 9 \Lambda (25 \Lambda - 16)$), and following \cite[Section~3]{DoubrovKruglikov}, choose $a$, $b$ so that $r = a^2 + b^2$, $s = a^2 b^2$. We fix a derivation $\mbF \in \Der(\mfn \otimes \bbC)$ (depending on $(r, s)$) that annihilates $\mfn_{-2} \otimes \bbC$ and preserves $\mfn_{-1} \otimes \bbC$. In particular, any such derivation (a) commutes with $\mbE$ and (b) can be identified with an element of the Lie algebra $\mfsp(\mfn_{-1} \otimes \bbC) \cong \mfsp(4, \bbC)$ of transformations preserving the canonical $(\mfn_{-2} \otimes \bbC)$-valued symplectic form $\Wedge^2 (\mfn_{-1} \otimes \bbC) \to \mfn_{-2} \otimes \bbC$ induced by the Lie bracket on $\mfn \otimes \bbC$. We then extend $\mfn \otimes \bbC$ by $\mbE, \mbF$ to produce the Lie algebra
 \begin{gather*}
	 \boxed{
		 \mfh_{r, s}^{\bbC} := (\mfn \otimes \bbC) \rightthreetimes \Span{\mbE, \mbF} \cong (\mfn \otimes \bbC) \rightthreetimes \bbC^2
 }
 \end{gather*}
occurring in the respective algebraic model. It turns out we are obliged to take $\mbF$ to correspond to an element of $\mfsp(4, \bbC)$ with eigenvalues $\pm a, \pm b$. This condition, together with the requirement that $\mfh_{r, s}^{\bbC}$ realizes an algebra model, turns out to determine $\mbF$ up to the natural action of $\Sp(4, \bbC)$. (The condition that $\mbD_{\Lambda}$ is not flat is equivalent to the requirement that $a \neq \pm 3b$ and $b \neq \pm 3a$.) We give for each model a Lie algebra isomorphism $\mfh_{r, s}^{\bbC} \cong \mfaut(\mbD_{r, s})$ identifying
 \begin{gather*}
 \mbE \leftrightarrow y \partial_y + p \partial_p + q \partial_q + 2 z \partial_z
 \qquad \textrm{and} \qquad
 \mbF \leftrightarrow \partial_x
 \end{gather*}
and identifying $\mbU$ with a constant multiple of $\partial_z$, and we set $\mfk := \Span{\mbE, \mbK}$ and $\mfd := \Span{\mbF, \mbL} \oplus \mfk$ for some complementary vectors $\mbK, \mbL$ depending on $(r, s)$.

It is convenient to treat separately the cases $\Lambda = 0$ (equivalently, $\mbF$ has a zero eigenvalue; we may take $r = 1$, $s = 0$, $a = 1$, $b = 0$) and $\Lambda = 1$ (equivalently, $\mbF$ has a repeated eigenvalue; we may take $r = 2$, $s = 1$, $a = b = 1$).

See Table \ref{table:complex-7} for details.

\subsection[Models with $\dim \mfh = 6$]{Models with $\boldsymbol{\dim \mfh = 6}$}
It remains to list the homogeneous complex distributions whose symmetry algebra has dimen\-sion~$6$. In addition to the data mentioned after Theorem~\ref{theorem:complex-classification}, Table~\ref{table:complex-6} also records for each algebraic model there a basis $(e_a)$ well-adapted to computing its admissible anti-involutions.
\begin{itemize}[leftmargin=1.5cm]\itemsep=0pt
\item[\textbf{$\mathbf{N.6}$}] In \cite{DoubrovGovorovCounterexample} Doubrov and Govorov reported a homogeneous distribution $\mbD$ with $\dim \mfaut(\mbD) = 6$ missing from Cartan's ostensible classification. Its symmetry algebra is a semidirect product $\mfsl_2(\bbC) \rightthreetimes (\mfm \otimes \bbC)$ of $\mfsl_2(\bbC)$ and the complex $3$-dimensional real Heisenberg algebra $\mfm \otimes \bbC$.
\item[\textbf{$\mathbf{D.6}_{\lambda}$}] These models, which are parameterized by $\lambda \in \bbC - \{ 0, \frac{1}{9}, 1, 9 \}$,\footnote{Taking $\lambda \in \{\frac{1}{9}, 9\}$ gives a flat distribution, and taking $\lambda \in \{0, 1\}$ yields a subspace $\mfd$ that does not satisfy the genericity criterion.} have symmetry algebra $\mfsl_2(\bbC) \oplus \mfsl_2(\bbC)$ \cite[Section~50]{CartanFiveVariables}. Where $\mbX, \mbY, \mbH$ and $\mbX', \mbY', \mbH'$ are standard bases of $\mfsl_2(\bbC)$ we define the models respectively by $\mfh := \mfsl_2(\bbC) \oplus \mfsl_2(\bbC)$, $\mfk: = \Span{\mbH + \mbH'}$, $\mfd := \Span{\mbX - \lambda \mbX', \mbY - \mbY'} \oplus \mfk$. For all parameter values $\lambda$, the Lie algebra automorphism of $\mfsl_2(\bbC) \oplus \mfsl_2(\bbC)$ interchanging the summands (interchanging $\mbX \leftrightarrow \mbX'$, $\mbY \leftrightarrow \mbY'$, $\mbH \leftrightarrow \mbH'$) is an isomorphism between $\mathbf{D.6}_{\lambda}$ and $\mathbf{D.6}_{1 / \lambda}$: It fixes $\mfk$ and $\Span{\mbY - \mbY'}$ and maps $\Span{\mbX - \lambda \mbX'}$ to $\Span{-\lambda \mbX + \mbX'} = \Span{\mbX - \lambda^{-1} \mbX'}$. On the other hand, Section~\ref{subsubsection:identification-6-complex} shows that parameter values $\lambda$, $\lambda'$ determine the same algebraic model iff $\lambda' = \lambda$ or $\lambda' = 1 / \lambda$.
\item[\textbf{$\mathbf{D.6}_{\infty}$}] This model has symmetry algebra $\mfsl_2(\bbC) \oplus \big(\mfso_2(\bbC) \rightthreetimes \bbC^2\big)$ \cite[Section~51]{CartanFiveVariables}.
\item[\textbf{$\mathbf{D.6}_{\ast}$}] This model has symmetry algebra $\mfso_3(\bbC) \rightthreetimes \bbC^3$ \cite[Section~52]{CartanFiveVariables}.
\end{itemize}

\begin{Example}\label{example:coordinates-to-distribution}
We indicate briefly by example how to product an algebraic model from a locally homogeneous distribution given in local coordinates.

Doubrov and Govorov gave model $\mathbf{N.6}$ in terms of the Monge normal form $F(x, y, p, q, z) := q^{1 / 3} + y$ (implicitly on an appropriate domain, say, on $\{q > 0\}$). From \cite{DoubrovGovorovCounterexample} the symmetry algebra $\mfh = \mfaut(\mbD_F)$ has basis
\begin{gather}
\begin{alignedat}{4}\label{equation:N.6-symmetries}
	\xi_1 &:= -y \partial_x + p^2 \partial_p + 3 p q \partial_q - \tfrac{1}{2} y^2 \partial_z , \qquad &
	\xi_4 &:= \partial_z , \\
	\xi_2 &:= -\big(x \partial_y + \partial_p + \tfrac{1}{2} x^2 \partial_z\big) , \qquad &
	\xi_5 &:= \partial_x , \\
	\xi_3 &:= -x \partial_x + y \partial_y + 2 p \partial_p + 3 q \partial_q , \qquad &
	\xi_6 &:= \partial_y + x \partial_z ,
\end{alignedat}
\end{gather}
but we can also compute $\mfh$ with the Maple package \texttt{DifferentialGeometry}:

\begin{Verbatim}[frame = single, xleftmargin=0.8cm, xrightmargin=0.8cm]
with(DifferentialGeometry): with(GroupActions):
DGsetup([x, y, p, q, z], M);

F := q^(1/3) + y;
Q := D_q;
X := evalDG(D_x + p * D_y + q * D_p + F * D_z);
DF := [Q, X];

InfinitesimalSymmetriesOfGeometricObjectFields([DF], output = "list");
\end{Verbatim}

The subalgebra $\Span{\xi_1, \xi_2, \xi_3}$ is isomorphic to $\mfsl_2(\bbC)$, and we may identify $(\xi_1, \xi_2, \xi_3)$ with the complexification of a standard basis $(\mbX, \mbY, \mbH)$ of $\mfsl_2$, namely one satisfying $[\mbX, \mbY] = \mbH, [\mbH, \mbX] = 2 \mbX, [\mbH, \mbY] = -2 \mbY$. The radical $\mfr$ of $\mfh$ is isomorphic to the complexification of the $3$-dimensional real Heisenberg algebra $\mfm$, and we may identify the basis $(\xi_4, \xi_5, \xi_6)$ of $\mfr$ with the complexification of a standard basis $(\mbU, \mbS, \mbT)$ thereof, namely one satisfying $[\mbS, \mbT] = \mbU$ and $[\mbU, \mbS] = [\mbU, \mbT] = 0$. Computing brackets realizes $\mfh$ as the complexification of the semidirect product $\mfsl_2 \rightthreetimes \mfm$ of $\mfsl_2$ and $\mfm$ specified by the bracket relations
\begin{gather*}
	\begin{array}{c|ccc}
		[\,\cdot\,,\,\cdot\,] & \mbU & \mbS &\mbT \\
		\hline
		\mbX & \cdot & \cdot & \mbS \\
		\mbY & \cdot & \mbT & \cdot \\
		\mbH & \cdot & \mbS & -\mbT
	\end{array} .
\end{gather*}

At the base point $u := (0, 0, 0, 1, 0) \in \bbC^5$, the isotropy subalgebra is $\mfk = \Span{\xi_1}$, and $\partial_q = \frac{1}{3} \xi_3$ and $D_x = -(\xi_2 - \xi_4 - \xi_5)$, so the resulting algebraic model $\mathbf{N.6}$ is specified by $\mfh = \mfsl_2(\bbC) \rightthreetimes (\mfm \otimes \bbC)$, $\mfk = \Span{\mbX}$, $\mfd = \Span{\mbY - \mbU - \mbS, \mbH} \oplus \mfk$.
\end{Example}

\section{Multiply transitive homogeneous real distributions}\label{section:classification-real}

We now use the list in Section~\ref{section:classification-complex} to classify multiply transitive homogeneous \textit{real} $(2, 3, 5)$ distributions.

\subsection{Real forms of complex algebraic models}

Given a real algebraic model $(\mfh, \mfk; \mfd)$, the triple $(\mfh \otimes \bbC, \mfk \otimes \bbC; \mfd \otimes \bbC)$ is a complex algebraic model; we call the latter the \textit{complexification} of the former.

Conversely, suppose that we have a complex algebraic model $(\mfh, \mfk; \mfd)$. Recall that a \textit{real form} of $\mfh$ is the fixed point Lie algebra $\mfh^{\phi}$ of an anti-involution $\phi \colon \mfh \to \mfh$, that is, a complex-antilinear map satisfying $\phi^2 = \id_{\mfh}$ and $\phi([x, y]) = [\phi(x), \phi(y)]$ for all $x, y, \in \mfh$. In analogy to \cite[Section~3]{DMT} we call $\phi$ \textit{admissible} iff it preserves $\mfk$ and $\mfd$, in which case $(\mfh^{\phi}, \mfk^{\phi}; \mfd^{\phi})$ is a real algebraic model, where $\mfk^{\phi} := \mfk \cap \mfh^{\phi}$ and $\mfd^{\phi} := \mfd \cap \mfh^{\phi}$. By construction its complexification is $(\mfh, \mfk; \mfd)$, so we call $(\mfh^{\phi}, \mfk^{\phi}; \mfd^{\phi})$ a \textit{real form} of $(\mfh, \mfk; \mfd)$.

We say that two admissible anti-involutions $\phi$, $\psi$ are \textit{equivalent} if $\psi = \alpha \circ \phi \circ \alpha^{-1}$ for some \textit{admissible} automorphism $\alpha$ of $\mfh$, that is, one preserving $\mfk$ and $\mfd$. Two admissible anti-involutions are equivalent iff they determine equivalent real algebraic models (and hence locally equivalent real homogeneous distributions), so to classify the latter one can classify the former. Not all $\mfh$ admit admissible anti-involutions, and hence not all complex algebraic models admit real forms.

We thus classify the multiply transitive homogeneous real distributions as follows: Up to local equivalence there is only a single real flat distribution and hence only a single real form of $\mbO$. For the submaximal case we appeal to the existing classification \cite[Theorem 2]{DoubrovZelenko} of real submaximal distributions.

This leaves the classification of real forms of the complex algebraic models $\dim \mfh = 6$. For each such model $(\mfh, \mfk; \mfd)$ (see Table \ref{table:complex-6}) we fix a basis $(e_1, \ldots, e_6)$ of $\mfh$ adapted to the filtration $\mfh \supset \mfd + [\mfd, \mfd] \supset \mfd \supset \mfk$ in the sense that
\begin{gather*}
	\mfk = \Span{e_6}, \qquad
	\mfd = \Span{e_4, e_5} \oplus \mfk, \qquad
	\mfd + [\mfd, \mfd] = \Span{e_3} \oplus \mfd ,
\end{gather*}
so that with respect to $(e_a)$ any admissible anti-involution $\phi$ of $\mfh$ is commensurately block lower-triangular; in particular, any admissible anti-involution $\phi$ satisfies $\phi(e_6) = \zeta e_6$ with $|\zeta| = 1$. Any such $\phi$ also preserves any other subspaces of $\mfh$ constructed invariantly from the data $(\mfh, \mfk; \mfd)$. In particular this includes the proper subspace $\mfe := [\mfk, \mfd] < \mfd$, but also in some cases centers and radicals, as well as brackets and intersections of other spaces constructed invariantly. In each case we are able to choose a basis well-adapted to some of these, and this restricts further in a~convenient way the form of $\phi$ with respect to the basis.

Next, for any automorphism $\alpha\colon \mfh \to \mfh$, by definition the constants $\sigma_{ab} := [\alpha(e_a), \alpha(e_b)] - \alpha([e_a, e_b])$ all vanish, and we impose those vanishing conditions to determine the admissible anti-involutions $\phi$ of the complex algebraic model. After classifying them up to equivalence, we record a representative anti-involution $\phi$, given in Tables \ref{table:real-N.6}--\ref{table:real-D.6-ast} with respect to the respective bases $(e_a)$, as well as the corresponding real model $(\mfh^{\phi}, \mfk^{\phi}; \mfd^{\phi})$. A reader interested in verifying the classification of admissible anti-involutions up to equivalence and the subsequent realization of the data given in those tables is encouraged to examine the Maple files accompanying 
this article.

\subsection{Local coordinate realizations}
For some purposes it is convenient to have local coordinate expressions of distributions. We give such forms for many of the real models in the classification, in some cases Monge normal forms, and indicate procedures for producing them in others.
\begin{enumerate}\itemsep=0pt
	\item For each complex algebraic model with $\dim \mfh = 6$ admitting a real form, the first representative anti-involution in the corresponding subsection fixes the adapted basis $(e_a)$ and hence its real span. So, if the data and Monge normal form specifying the complex algebraic model can be interpreted as real, doing so gives a real algebraic model and a corresponding real Monge normal form. This procedure immediately yields Monge normal forms for models $\mathbf{N.6}^+$, $\mathbf{D.6}_\lambda^{2-}$ ($\lambda > 0$), $\mathbf{D.6}_\infty^1$, and $\mathbf{D.6}_\ast^{1-}$. The conclusion applies just as well to the real forms $\mbO^{\bbR}$, $\mathbf{N.7}_0^{\rightarrow}$, and $\mathbf{N.7}_\Lambda^\nearrow$, $\Lambda \not\in [0, 1)$.
	\item The submaximal real distributions $\mathbf{N.7}_\Lambda^\bullet$ were classified in \cite{DoubrovKruglikov}, and Monge normal forms were recorded there. See Section~\ref{subsection:N.7-real}.
	\item Example \ref{example:D.6-ast-coordinates} later in this section outlines, using the real model $\mathbf{D.6}_{\ast}^4$ as an example, how to construct local coordinates (and indeed, a Monge normal form) from an algebraic model.
	\item Example \ref{example:N.6-Monge} applies the identification algorithm in that section to show that a particular function $F$ defines a Monge normal form for the real algebraic model $\mathbf{N.6}^-$.
	\item Section \ref{section:rolling} realizes several of the real models in the classification as rolling distributions, from which one can readily construct coordinate realizations; this is carried out explicitly for the real algebraic models $\mathbf{D.6}_{\lambda}^6$.
\end{enumerate}


\subsection{The real flat model \texorpdfstring{$\mbO^{\bbR}$}{O-R}}\label{subsection:O-real}

As in the complex case, up to local equivalence there is a unique flat distribution. Thus all admissible anti-involutions of $\mfg_2(\bbC)$ are equivalent; taking complex conjugation with respect to the realization \eqref{equation:g2} yields the model
\begin{gather*}
	\boxed{
		\big(\mfg_2, \mfq; \mfg_2^{-1}\big) ,
	}
\end{gather*}
which we denote $\mbO^{\bbR}$.

The real flat model can be described as the rolling distribution (see Section~\ref{section:rolling}) determined by a pair of spheres, one whose radius is thrice that of the other \cite{BorMontgomery},\footnote{Agrachev \cite[Section~1]{Agrachev} attributes this characterization to Bryant, who pointed out in a note to Bor and Montgomery excerpted in the introduction of \cite{BorMontgomery} that it can be derived from a characterization due to Cartan \cite[Section~53]{CartanFiveVariables}.} and also as a canonical distribution determined by the algebra $\wt{\bbO}$ of split octonions on the null quadric in the projectivization $\bbP\big(\operatorname{Im} \wt{\bbO}\big)$ of the space of purely imaginary split octonions \cite[Section~6]{BorMontgomery}, \cite{Sagerschnig}; see also \cite{BaezHuerta}.

\subsection{The submaximal real models \texorpdfstring{$\mathbf{N.7}_\Lambda^\bullet$}{N.7-Lambda}}\label{subsection:N.7-real}

The classification of real submaximally symmetric $(2, 3, 5)$ distributions was established in \cite[Theorem 2]{DoubrovZelenko} using the geometry of distinguished curves called abnormal extremals: Any such distribution can be written in the Monge normal form
\begin{gather*}
	F_{r, s}(x, y, p, q, z) = q^2 + r p^2 + s y^2
\end{gather*}
for some $r, s \in \bbR$, and as in the complex case, the distribution $\mbD_{r, s}$ so determined is locally equivalent to the flat model iff $9 r^2 = 100 s$, which again we henceforth exclude. Otherwise, $\mbD_{r, s}$ and $\mbD_{r', s'}$ are locally equivalent iff there is a constant $c \in \bbR_+$ such that $r' = c r$, $s' = c^2 s$.

If the complex algebraic model $\mathbf{N.7}_\Lambda$ admits a real form, then it follows from the algorithm in Section~\ref{subsubsection:identification-submaximal-complex} that $\Lambda$ is real; conversely, reality of $\Lambda$ turns out also to be a sufficient condition for the existence of a real form. The form of the equivalence relation implies that the triple $(\Lambda, \sign(r), \sign(s))$, is a complete invariant of the model. We decorate our model names~$\mathbf{N.7}_\Lambda^\bullet$ by replacing $\bullet$ with an arrow (one of $\rightarrow$, $\nearrow$, $\uparrow$, $\nwarrow$, $\leftarrow$, $\swarrow$, $\downarrow$, $\searrow$) to indicate the pair $(\sign(r), \sign(s))$.

Recall that the explicit complex algebraic models $\mathbf{N.7}_\Lambda$ in Section~\ref{subsection:N.7} were specified using constants $a$, $b$ satisfying $r = a^2 + b^2$, $s = a^2 b^2$ for a pair $(r, s)$ from which we could recover $\Lambda$. This has the advantage that we can give explicit algebraic models for all submaximal models with little case splitting, but for $r, s \in \bbR$, the corresponding parameters $a$, $b$ are real iff $r, s \geq 0$, $r^2 \geq 4 s$. Thus, we record explicit data defining the real algebraic models $\mathbf{N.7}_\Lambda$ in Table \ref{table:real-N.7} in Appendix~\ref{appendixB}. The elements $\mbF$ defined there are elements of $\Der(\mfn)$.

\subsection[The real forms of model $\mathbf{N.6}$]{The real forms of model $\boldsymbol{\mathbf{N.6}}$}\label{subsection:N.6-real}
The adapted basis $(e_a)$ is
\begin{alignat*}{5}
&	e_1 := 3 \mbT , \qquad && e_2 := \mbU , \qquad && 	e_3 := 2 \mbU + 3 \mbS , \qquad && 	e_4 := \mbY - \mbU - \mbS , & \\
&	e_5 := \mbH , \qquad && 	e_6 := \mbX . && && &
\end{alignat*}
The center of $\mfh$ is $\Span{\mbU} = \Span{e_2}$, so $\phi(e_2) = \beta e_2$ for some $\beta$. The radical $\mfr$ of $\mfh$ is $\Span{\mbS, \mbT, \mbU} = \Span{e_1, e_2, e_3} \cong \mfm \otimes \bbC$, so $\phi(e_1) = \alpha_1 e_1 + \alpha_2 e_2 + \alpha_3 e_3$, and $\mfr \cap (\mfd + [\mfd, \mfd]) = \Span{e_3}$, so $\phi(e_3) = \gamma e_3$. Finally, $\mfe := [\mfk, \mfd] = \Span{\mbX, \mbH} = \Span{e_5, e_6}$, so $\phi(e_5) = \epsilon_5 e_5 + \epsilon_6 e_6$, and adaptation of the basis gives that $\phi(e_4) = \delta_4 e_4 + \delta_5 e_5 + \delta_6 e_6$. Nondegeneracy implies that $\alpha_1, \beta, \gamma, \delta_4, \epsilon_5 \neq 0$.

Now, $\sigma_{56} = 0$ implies that $\epsilon_5 = 1$, and then $\sigma_{35} = \sigma_{45} = 0$ implies $\beta = \gamma = \delta_4$. Next, $\sigma_{13} = \sigma_{16} = \sigma_{34} = 0$ implies that $\alpha_1 = 1$ and $\delta_4 = \zeta = \pm 1$, then $\sigma_{15} = 0$ implies $\alpha_2 = -2 \alpha_3 = \epsilon_6$, and then $\sigma_{14} = 0$ implies $\delta_6 = \mp \delta_5^2$ and $\epsilon_6 = \mp 2 \delta_5$; in summary,
\begin{alignat*}{3}
	 \phi(e_1) &= e_1 \mp 2 \delta e_2 \pm \delta e_3, \qquad
	&\phi(e_2) &= \pm e_2, \qquad
	&\phi(e_3) &= \pm e_3, \\
	 \phi(e_4) &= \pm e_4 + \delta e_5 \mp \delta^2 e_6, \qquad
	&\phi(e_5) &= e_5 \mp 2 \delta e_6, \qquad
	&\phi(e_6) &= \pm e_6 ,
\end{alignat*}
where $\delta := \delta_5$. Conjugating one automorphism of this form by another preserves the sign $\pm$, so admissible anti-involutions with different choices of sign are inequivalent.

\smallskip
\noindent \textbf{Case: $\zeta = + 1$.} Computing directly gives that for $\phi$ complex-antilinear the condition $\phi^2 = \id$ is equivalent to $\Re \delta = 0$, and conjugating $\phi$ by the admissible automorphism $(e_1, \ldots, e_6) \mapsto \big(e_1 - 2 t e_2 + t e_3, e_2, e_3, e_4 + t e_5 - t^2 e_6, e_5 - 2 t e_6, e_6\big)$ induces the transformation $\delta \rightsquigarrow \delta + 2 {\rm i} \Im t$. Setting $t = -\frac{1}{2} \delta$ shows that all admissible anti-involutions in this branch are equivalent to the one with $\delta = 0$, for which the fixed point Lie algebra $\mfh^{\phi} = \spanOperator_{\bbR}\{e_1, \ldots, e_6\}$. The corresponding real form, which we denote $\mathbf{N.6}^+$, is
\begin{gather*}
	\boxed{
		\mfh^{\phi} = \mfsl_2 \rightthreetimes \mfm , \qquad
		\mfk^{\phi} = \Span{\mbX}, \qquad
		\mfd^{\phi} = \Span{\mbY - \mbU - \mbS, \mbH} \oplus \mfk .
	}
\end{gather*}

\noindent \textbf{Case: $\zeta = - 1$.} In this case $\phi$ is an anti-involution iff $\Im \delta = 0$, and conjugating by the admissible automorphism $(e_1, \ldots, e_6) \mapsto \big(e_1 - 2 t e_2 + t e_3, -e_2, -e_3, -e_4 - t e_5 + t^2 e_6, e_5 - 2 t e_6, -e_6\big)$ induces $\delta \rightsquigarrow \delta + 2 \Re t$, and so we may normalize to $\delta = 0$. For this $\phi$, $\mfh^{\phi} = \spanOperator_{\bbR}\{e_1, {\rm i} e_2, {\rm i} e_3, {\rm i} e_4, e_5, {\rm i} e_6\}$, and we may identify $\mfh^{\phi} = \mfsl_2 \rightthreetimes \mfm$
via
\begin{gather*}
	\mbX \leftrightarrow {\rm i} e_6 , \qquad
	\mbY \leftrightarrow -\tfrac{1}{3} {\rm i} (e_2 + e_3 + 3 e_4) , \qquad
	\mbH \leftrightarrow e_5 , \qquad
	\mbU \leftrightarrow {\rm i} e_2 , \\
	\mbS \leftrightarrow \tfrac{1}{3} {\rm i} (-2 e_2 + e_3) , \qquad
	\mbT \leftrightarrow \tfrac{1}{3} e_1 ,
\end{gather*}
where the semidirect product is as in the case $\zeta = + 1$. The model, which we denote $\mathbf{N.6}^-$, is
\begin{gather*}
	\boxed{
 \mfh^{\phi} = \mfsl_2 \rightthreetimes \mfm , \qquad
		\mfk^{\phi} = \Span{\mbX} , \qquad
		\mfd^{\phi} = \Span{\mbY + \mbU + \mbS, \mbH} \oplus \mfk^{\phi} .
	}
\end{gather*}
In Example \ref{example:N.6-Monge} we show that $F(x, y, p, q, z) = q^{1 / 3} - y$ realizes $\mathbf{N.6}^-$ in Monge normal form.

\subsection{The real forms of models \texorpdfstring{$\mathbf{D.6}_{\lambda}$}{D.6-lambda}}\label{subsection:D.6-lambda-real}

This case is the most involved. We will see that (1) in some cases the qualitative features of the real forms (including the isomorphism types of the real forms $\mfh^{\phi}$) depend on the sign of $\lambda$, and (2) the case $\lambda = -1$ is distinguished.

In Section~\ref{subsubsection:identification-6-complex} we show that $\lambda$ is an invariant of the distribution up to inversion, and it is manifestly real for distributions that are complexifications of real distributions. Thus, only the models $\mathbf{D.6}_{\lambda}$ with $\lambda$ real can admit real forms, and we therefore restrict to such $\lambda$. Since for real $\lambda$ the abstract data defining the model $\mathbf{D.6}_{\lambda}$ can be interpreted as real, all such models do admit real forms.

Before proceeding, recall the classification of the real forms of $\mfsl_2(\bbC) \oplus \mfsl_2(\bbC) \cong \mfso_4(\bbC)$ and the fact that those real forms are determined up to isomorphism by the signatures of their Killing forms. We will see that all four forms occur in the classification.\footnote{Here $\mfsl_2(\bbC)_{\bbR}$ is the real Lie algebra underlying $\mfsl_2(\bbC)$.}

\begin{table}[htb]\small\centering
\caption{Real forms of $\mfsl_2(\bbC) \oplus \mfsl_2(\bbC) \cong \mfso_4(\bbC)$.}\label{table:real-forms-sl2C-sl2C}

\vspace{1mm}

\begin{tabular}{cc}
\toprule
 real form
 & signature \\
\midrule
	 $\mfso_3 \oplus \mfso_3 \cong \mfso_4$
 & $(0, 6)$ \\
	 $\mfsl_2 \oplus \mfso_3 \cong \mfso_4^*$
 & $(2, 4)$ \\
	 $\mfsl_2(\bbC)_{\bbR} \cong \mfso_{1, 3}$
 & $(3, 3)$ \\
	 $\mfsl_2 \oplus \mfsl_2 \cong \mfso_{2, 2}$
 & $(4, 2)$ \\
\bottomrule
\end{tabular}
\end{table}

The adapted basis $(e_a)$,
\begin{equation*}
\begin{alignedat}{3}
	 e_1 &:= \mbX - \lambda^2 \mbX' , \qquad
	&e_2 &:= \mbY - \lambda \mbY' , \qquad
	&e_3 &:= \mbH + \lambda \mbH' , \\
	 e_4 &:= \mbX - \lambda \mbX' , \qquad
	&e_5 &:= \mbY - \mbY' , \qquad
	&e_6 &:= \mbH + \mbH' ,
\end{alignedat}
\end{equation*}
satisfies $\mfe := [\mfk, \mfd] = \Span{e_4, e_5}$, $[\mfe, \mfe] = \Span{e_3}$, and $[\mfe, [\mfe, \mfe]] = \Span{e_1, e_2}$. Since $\phi$ preserves each of these subspaces, it has the form
\begin{gather}\label{equation:phi-D.6-lambda}
\begin{alignedat}{3}
	 \phi(e_1) &= \alpha_1 e_1 + \alpha_2 e_2, \qquad
	&\phi(e_2) &= \beta_1 e_1 + \beta_2 e_2, \qquad
	&\phi(e_3) &= \gamma e_3 , \\
	 \phi(e_4) &= \delta_4 e_4 + \delta_5 e_5, \qquad
	&\phi(e_5) &= \epsilon_4 e_4 + \epsilon_5 e_5, \qquad
	&\phi(e_6) &= \zeta e_6 ,
\end{alignedat}
\end{gather}
for some coefficients, and nondegeneracy implies $\alpha_1 \beta_2 - \alpha_2 \beta_1, \gamma, \delta_4 \epsilon_5 - \delta_5 \epsilon_4 \neq 0$.

Now, $\sigma_{34} = \sigma_{35} = 0$ implies that $\alpha_1 = \gamma \delta_4$, $\alpha_2 = -\gamma \delta_5$, $\beta_1 = -\gamma \epsilon_4$, $\beta_2 = \gamma \epsilon_5$, then $\sigma_{16} = 0$ forces $\zeta = \pm' 1$, and then $\sigma_{13} = 0$ implies $\gamma = \pm 1$. Conjugating one automorphism of this form by another shows that the respective signs $\pm$, $\pm'$ of $\gamma$, $\zeta$ are invariant under this conjugation, so admissible anti-involutions with different choices of $\pm$, $\pm'$ are inequivalent.

Forming a suitable linear combination of the $e_3$ and $e_6$ components of $\sigma_{15}$ gives that $(\lambda + 1)(\gamma - \zeta) = 0$, so $\gamma = \zeta$ or $\lambda = -1$. For each subcase, the conditions $\sigma_{16} = \sigma_{26} = 0$ imply the vanishing either of both $\delta_5$ and $\epsilon_4$ or of both $\delta_4$ and $\epsilon_5$, then $\sigma_{12} = 0$ implies that one of the two remaining quantities can be written in terms of the other, after which all of the conditions $\sigma_{ab} = 0$ are satisfied.

\noindent \textbf{Case: $\gamma = \zeta$.} We split cases according to the sign of $\zeta = \pm 1$.
\begin{itemize}[leftmargin=\caseIndent]\itemsep=0pt
	\item[] \textbf{Subcase: $\zeta = +1$.} We have $\delta_5 = \epsilon_4 = 0$ and $\epsilon_5 = \delta_4^{-1}$, so that $\phi(e_1) = \delta e_1$, $\phi(e_2) = \delta^{-1} e_2$, $\phi(e_3) = e_3$, $\phi(e_4) = \delta e_4$, $\phi(e_5) = \delta^{-1} e_5$, $\phi(e_6) = e_6$, where $\delta := \delta_4$, and the condition that $\phi$ is an anti-involution is $|\delta| = 1$. The admissible automorphism $(e_1, \ldots, e_6) \mapsto \big(t e_1, t^{-1} e_2, e_3, t e_4, e^{-1} e_5, e_6\big)$ induces $\delta \rightsquigarrow t^2 \delta / |t|^2$, so we may normalize to $\delta = 1$, for which $\mfh^{\phi} = \spanOperator_{\bbR}\{e_1, \ldots, e_6\}$. The gives a model, $\mathbf{D.6}_{\lambda}^{2-}$:
	\begin{gather*}
		\boxed{
			\mfh^{\phi} = \mfsl_2 \oplus \mfsl_2, \qquad
			\mfk^{\phi} = \Span{\mbH + \mbH'} , \qquad
			\mfd^{\phi} = \Span{\mbX - \lambda \mbX', \mbY - \mbY'} \oplus \mfk^{\phi} .
		}
	\end{gather*}
	As in the complex case, the Lie algebra isomorphism exchanging the direct summands $\mfsl_2$ defines an isomorphism between $\mathbf{D.6}_{\lambda}^{2-}$ and $\mathbf{D.6}_{1 / \lambda}^{2-}$; complexifying shows that this again exhausts the isomorphisms among these models.
	\item[] \textbf{Subcase: $\zeta = -1$.} Now, $\delta_4 = \epsilon_5 = 0$, $\epsilon_4 = \delta_5^{-1}$, so $\phi(e_1) = \delta e_2$, $\phi(e_2) = \delta^{-1} e_1$, $\phi(e_3) = -e_3$, $\phi(e_4) = \delta e_5$, $\phi(e_5) = \delta^{-1} e_4$, $\phi(e_6) = -e_6$, where $\delta := \delta_5$, and the anti-involution condition is $\Im \delta = 0$. The admissible automorphism $(e_1, \ldots, e_6) \mapsto \big(t e_2, t^{-1} e_1, -e_3, t e_5, t^{-1} e_4, -e_6\big)$ induces $\delta \rightsquigarrow |t|^2 / \delta$, so we may normalize to $\delta = \pm 1$. The two anti-involutions these values determine give rise to real forms $\mfh^{\phi}$ whose Killing forms have different signatures, so they cannot be equivalent.
\begin{itemize}[leftmargin=\caseIndent]\itemsep=0pt
	\item[] \textbf{Subsubcase: $\delta = +1$.} In this case $\mfh^{\phi} = \spanOperator_{\bbR}\{e_1 + e_2, {\rm i} (e_1 - e_2), {\rm i} e_3, e_4 + e_5, {\rm i} (e_4 - e_5), {\rm i} e_6\}$.
		\begin{itemize}[leftmargin=\caseIndent]\itemsep=0pt
		\item[] \textbf{Subsubsubcase: $\lambda > 0$.} The signature of the Killing form of $\mfh^{\phi}$ is $(4, 2)$, so $\mfh^{\phi} \cong \mfsl_2 \oplus \mfsl_2$, and we can realize this isomorphism via
			\begin{align*}
				\mbX\phantom{'}
					&\leftrightarrow \tfrac{1}{\sqrt{2} (\lambda - 1)} (e_1 + e_2 + {\rm i} e_3 - \lambda e_4 - \lambda e_5 - {\rm i} \lambda e_6) , \\
				\mbY\phantom{'}
					&\leftrightarrow \tfrac{1}{\sqrt{2} (\lambda - 1)} {\rm i} (-e_1 + e_2 - e_3 + \lambda e_4 - \lambda e_5 + \lambda e_6), \\
				\mbH\phantom{'}
					&\leftrightarrow \tfrac{1}{\lambda - 1} [(-1 - {\rm i}) e_ 1 + (-1 + {\rm i}) e_2 - {\rm i} e_3 + (1 + {\rm i}) \lambda e_4 + (1 - {\rm i}) \lambda e_5 + {\rm i} \lambda e_6], \\
 \mbX'
					&\leftrightarrow \tfrac{1}{(\lambda - 1) \sqrt{2 \lambda}} \big({-}e_1 - e_2 - {\rm i} \sqrt{\lambda} e_3 + e_4 + e_5 + {\rm i} \sqrt{\lambda} e_6\big), \\
				\mbY'
					&\leftrightarrow \tfrac{1}{(\lambda - 1) \sqrt{2 \lambda}} {\rm i} \big(e_1 - e_2 + \sqrt{\lambda} e_3 - e_4 + e_5 - \sqrt{\lambda} e_6\big), \\
				\mbH'
					&\leftrightarrow \tfrac{1}{(\lambda - 1) \sqrt{\lambda}} \big[(1 + {\rm i}) e_1 + (1 - {\rm i}) e_2 + {\rm i} \sqrt{\lambda} e_3 \\
 &\quad{} + (-1 - {\rm i}) e_4 + (-1 + {\rm i}) e_5 - \sqrt{\lambda} e_6\big] ,
			\end{align*}
		and then the model, which we denote $\mathbf{D.6}_{\lambda}^{2+}$, is
			\begin{equation*}
				\boxed{
 \begin{split}
& \mfh^{\phi} = \mfsl_2 \oplus \mfsl_2 , \qquad
	 	 \mfk^{\phi} = \operatorname{span}\big\{\wt\mbH + \wt\mbH'\big\} , \\
& \mfd^{\phi} = \operatorname{span}\big\{\wt\mbX + \sqrt{\lambda} \wt\mbX', \wt\mbY + \sqrt{\lambda} \wt\mbY'\big\} \oplus \mfk^{\phi} ,
 \end{split}
				}
			\end{equation*}
		where\footnote{The basis $(\wt\mbX, \wt\mbY, \wt\mbH)$ of $\mfsl_2 \cong \mfso_{1, 2}$ is pseudo-orthonormal with respect to an appropriate multiple of the Killing form.}
		\begin{gather*}
			\wt\mbX := \tfrac{1}{\sqrt{2}} \mbX + \tfrac{1}{2} \mbH , \qquad
			\wt\mbY := -\tfrac{1}{\sqrt{2}} \mbY + \tfrac{1}{2} \mbH , \qquad
			\wt\mbH := \tfrac{1}{\sqrt{2}} (\mbX - \mbY) + \tfrac{1}{2} \mbH ,
		\end{gather*}
		and $\wt\mbX'$, $\wt\mbY'$, $\wt\mbH'$ are defined analogously. Exchanging the direct summands $\mfsl_2$ is an isomorphism between $\mathbf{D.6}_{\lambda}^{2+}$ and $\mathbf{D.6}_{1 / \lambda}^{2+}$, and this exhausts the isomorphisms among these models.

		\item[] \textbf{Subsubsubcase: $\lambda < 0$.} The Killing form has signature $(2, 4)$, so $\mfh^{\phi} \cong \mfsl_2 \oplus \mfso_3$. We can identify $\mbX, \mbY, \mbH \in \mfsl_2$ as in the case $\lambda > 0$, and if we denote by $(\mbA, \mbB, \mbC)$ a standard basis of $\mfso_3(\bbR)$~-- one characterized by $[\mbA, \mbB] = \mbC$ and its cyclic permutations~-- we can complete the identification via
		\begin{align*}
			\mbA &\leftrightarrow \tfrac{1}{2 (\lambda - 1) \sqrt{-\lambda}} {\rm i} (e_1 - e_2 - e_4 + e_5) , \\
			\mbB &\leftrightarrow \tfrac{1}{2 (\lambda - 1) \sqrt{-\lambda}} (e_1 + e_2 - e_4 - e_5) , \\
			\mbC &\leftrightarrow \tfrac{1}{2 (\lambda - 1)} {\rm i} (-e_3 + e_6) ,
		\end{align*}
		so
			\begin{equation*}
				\boxed{
 \begin{split}
 & \mfh^{\phi} = \mfsl_2 \oplus \mfso_3 , \qquad
	 	 \mfk^{\phi} = \operatorname{span}\big\{\wt\mbH + \mbC\big\} , \\
& \mfd^{\phi} = \operatorname{span}\big\{
								\wt\mbX + \sqrt{-\lambda} \mbA ,
								\wt\mbY + \sqrt{-\lambda} \mbB
							 \big\} \oplus \mfk^{\phi} ,
 \end{split}
				}
			\end{equation*}
		We denote this model $\mathbf{D.6}_{\lambda}^4$.
		\end{itemize}

	\item[] \textbf{Subsubcase: $\delta = -1$.} In this case $\mfh^{\phi} = \spanOperator_{\bbR}\{ e_1 - e_2, {\rm i} (e_1 + e_2), {\rm i} e_3, e_4 - e_5, {\rm i} (e_4 + e_5), {\rm i} e_6 \}$.
		\begin{itemize}[leftmargin=\caseIndent]\itemsep=0pt
			\item[] \textbf{Subsubsubcase: $\lambda > 0$.} The Killing form of $\mfh^{\phi}$ is definite, so $\mfh^{\phi} \cong \mfso_3 \oplus \mfso_3$, and we can realize this isomorphism via
		\begin{alignat*}{4}
			\mbA
				& \leftrightarrow \tfrac{1}{2 (\lambda - 1)} (e_1 - e_2 - \lambda e_4 + \lambda e_5) , \qquad &
			\mbA'
				&\leftrightarrow \tfrac{1}{2 (\lambda - 1) \sqrt{\lambda}} (e_1 - e_2 - e_4 + e_5) , \\
			\mbB
				&\leftrightarrow \tfrac{1}{2 (\lambda - 1)} {\rm i} (e_1 + e_2 - \lambda e_4 - \lambda e_5) , \qquad &
			\mbB'
				&\leftrightarrow \tfrac{1}{2 (\lambda - 1) \sqrt{\lambda}} {\rm i} (e_1 + e_2 - e_4 - e_5) , \\
			\mbC
				&\leftrightarrow \tfrac{1}{2 (\lambda - 1)} {\rm i} (-e_3 + \lambda e_6 ) , \qquad &
			\mbC'
				&\leftrightarrow \tfrac{1}{2 (\lambda - 1)} {\rm i} (e_3 - e_6) ,
		\end{alignat*}
		so the model, which we denote $\mathbf{D.6}_{\lambda}^6$, is
		\begin{equation*}
			\boxed{
 \begin{split}
 &\mfh^{\phi} = \mfso_3 \oplus \mfso_3 , \qquad
	 	 \mfk^{\phi} = \operatorname{span}\big\{\mbC + \mbC'\big\} , \\
 &\mfd^{\phi} = \operatorname{span}\big\{\mbA - \sqrt{\lambda} \mbA', \mbB - \sqrt{\lambda} \mbB'\big\} \oplus \mfk^{\phi} ,
 \end{split}
			}
		\end{equation*}
		where $(\mbA', \mbB', \mbC')$ is a standard basis of the second summand $\mfso_3$. Exchanging the direct summands $\mfso_3$ is an isomorphism between $\mathbf{D.6}_{\lambda}^6$ and $\mathbf{D.6}_{1 / \lambda}^6$, and this exhausts the isomorphisms among these models. Example~\ref{example:two-spheres} realizes the models $\mathbf{D.6}_{\lambda}^6$ in local coordinates.

		\item[] \textbf{Subsubsubcase: $\lambda < 0$.} The Killing form of $\mfh^{\phi}$ has signature $(2, 4)$, so $\mfh^{\phi} \cong \mfso_3 \oplus \mfsl_2$, and we can realize this isomorphism by identifying $\mbA$, $\mbB$, $\mbC$ as in the $\lambda > 0$ case and identifying
		\begin{align*}
			\mbX
				&\leftrightarrow \tfrac{1}{(\lambda - 1) \sqrt{- 2 \lambda}} {\rm i} \big({-}e_1 - e_2 - \sqrt{-\lambda} e_3 + e_4 + e_5 + {\rm i} \sqrt{-\lambda} e_6\big) , \\
			\mbY
				&\leftrightarrow \tfrac{1}{(\lambda - 1) \sqrt{- 2 \lambda}} \big(e_1 - e_2 - {\rm i} \sqrt{-\lambda} e_3 - e_4 + e_5 + {\rm i} \sqrt{-\lambda} e_6\big) , \\
			\mbH
				&\leftrightarrow \tfrac{1}{(\lambda - 1) \sqrt{-\lambda}} \big[(1 + {\rm i}) e_1 + (-1 + {\rm i}) e_2 - {\rm i} \sqrt{-\lambda} e_3 \\
 &\quad {}+ (-1 - {\rm i}) e_4 + (1 - {\rm i}) e_5 + {\rm i} \sqrt{-\lambda} e_6\big] .
		\end{align*}
		Then,
		\begin{equation*}
			\boxed{
 \begin{split}
 &\mfh^{\phi} = \mfso_3 \oplus \mfsl_2 , \qquad
	 	 \mfk^{\phi} = \operatorname{span}\big\{\mbC + \wt\mbH\big\} , \\
 &\mfd^{\phi} = \operatorname{span}\big\{
								\mbA + \sqrt{-\lambda} \wt\mbX,
								\mbB + \sqrt{-\lambda} \wt\mbY
							 \big\} \oplus \mfk^{\phi} ,
 \end{split}
			}
		\end{equation*}
		The isomorphism $\mfsl_2 \oplus \mfso_3 \cong \mfso_3 \oplus \mfsl_2$ given by reversing the order of the factors identifies the model with parameter $\lambda$ with $\mathbf{D.6}_{1 / \lambda}^4$, so this branch yields no new models.
		\end{itemize}
	\end{itemize}
\end{itemize}

\smallskip
\noindent \textbf{Case: $\gamma \neq \zeta$.} We have $-\gamma = \zeta = \pm 1$.
\begin{itemize}[leftmargin=\caseIndent]
	\item[] \textbf{Subcase: $\zeta = +1$.} Here, $\delta_5 = \epsilon_4 = 0$, $\epsilon_5 = -\delta_4^{-1}$, so $\phi(e_1) = -\delta e_1$, $\phi(e_2) = \delta^{-1} e_2$, $\phi(e_3) = -e_3$, $\phi(e_4) = \delta e_4$, $\phi(e_5) = -\delta^{-1} e_5$, $\phi(e_6) = e_6$, where $\delta := \delta_4$, and the anti-involution condition is $|\delta| = 1$. The admissible automorphism $(e_1, \ldots, e_6) \mapsto \big({-}t e_1, t^{-1} e_2, e_3, t e_4, -t^{-1} e_5, e_6\big)$ induces $\delta \rightsquigarrow t^2 \delta / |t|^2$; normalizing to $\delta = 1$ gives $\mfh^{\phi} = \spanOperator_{\bbR}\{{\rm i} e_1, e_2, {\rm i} e_3, e_4, {\rm i} e_5, e_6\}$. The Killing form of $\mfh^{\phi}$ has signature $(3, 3)$, so $\mfh^{\phi} \cong \mfso_{1, 3}$.

 If we fix a basis $(\mbA, \mbB, \mbC, \mbD_{\mbA}, \mbD_{\mbB}, \mbD_{\mbC})$ of $\mfso_{1, 3}$ satisfying the identities $[\mbA, \mbB] = \mbC$, $[\mbA, \mbD_{\mbB}] = -[\mbB, \mbD_{\mbA}] = \mbD_{\mbC}$, and $[\mbD_{\mbA}, \mbD_{\mbB}] = \mbC$ and their cyclic permutations, as well as $[\mbA, \mbD_{\mbA}] = [\mbB, \mbD_{\mbB}] = [\mbC, \mbD_{\mbC}] = 0$, then we may realize this identification explicitly via
		\begin{alignat*}{4}
			 \mbA &\leftrightarrow \tfrac{1}{2} (- e_2 + e_4) , \qquad
			&\mbD_{\mbA} &\leftrightarrow \tfrac{1}{2} {\rm i} (- e_1 + e_5) , \\
			 \mbB &\leftrightarrow \tfrac{1}{2} {\rm i} ( e_1 + e_5) , \qquad
			&\mbD_{\mbB} &\leftrightarrow \tfrac{1}{2} ( e_2 + e_4) , \\
			 \mbC &\leftrightarrow \tfrac{1}{2} {\rm i} e_3 , \qquad
			&\mbD_{\mbC} &\leftrightarrow \tfrac{1}{2} e_6 .
		\end{alignat*}
		Then, the model, which we denote $\mathbf{D.6}_{-1}^{3-}$, is
		\begin{gather*}
			\boxed{
 \mfh^{\phi} \cong \mfso_{1, 3} , \qquad
 \mfk^{\phi} = \Span{\mbD_{\mbC}} , \qquad
				\mfd^{\phi} = \Span{\mbA + \mbD_{\mbB}, \mbB + \mbD_{\mbA}} \oplus \mfk^{\phi} .
			}
		\end{gather*}
	\item[] \textbf{Subcase: $\zeta = -1$.} Here, $\delta_4 = \epsilon_5 = 0$ and $\epsilon_4 = -\delta_5^{-1}$, so $\phi(e_1) = -\delta e_2$, $\phi(e_2) = \delta^{-1} e_1$, $\phi(e_3) = e_3$, $\phi(e_4) = \delta e_5$, $\phi(e_5) = -\delta^{-1} e_4$, $\phi(e_6) = -e_6$, where $\delta := \delta_5$, and the anti-involution condition is $\Re \delta = 0$. The admissible automorphism $(e_1, \ldots, e_6) \mapsto \big({-}t e_2, t^{-1} e_1, e_3, t e_5, -t^{-1} e_4, -e_6\big)$ induces $\delta \rightsquigarrow |t|^2 / \delta$, so we may normalize to $\delta = {\rm i}$, for which $\mfh^{\phi} = \spanOperator_{\bbR}\{e_1 - {\rm i} e_2, {\rm i} e_1 - e_2, e_3, e_4 + {\rm i} e_5, {\rm i} e_4 + e_5, {\rm i} e_6\}$. The Killing form of $\mfh^{\phi}$ has signature $(3, 3)$, so again $\mfh^{\phi} \cong \mfso_{1, 3}$, and we may realize this identification via
		\begin{alignat*}{4}
			 \mbA &\leftrightarrow \tfrac{1}{2 \sqrt{2}} ({\rm i} e_1 - e_2 + e_4 + {\rm i} e_5) , \qquad
			&\mbD_{\mbA} &\leftrightarrow \tfrac{1}{2 \sqrt{2}} (-{\rm i} e_1 + e_2 + e_4 + {\rm i} e_5) , \\
			 \mbB &\leftrightarrow \tfrac{1}{2 \sqrt{2}} (- e_1 + {\rm i} e_2 + {\rm i} e_4 + e_5) , \qquad
			&\mbD_{\mbB} &\leftrightarrow \tfrac{1}{2 \sqrt{2}} (e_1 - {\rm i} e_2 + {\rm i} e_4 + e_5) , \\
			 \mbC &\leftrightarrow \tfrac{1}{2} {\rm i} e_6 , \qquad
			&\mbD_{\mbC} &\leftrightarrow \tfrac{1}{2} e_3 .
		\end{alignat*}
		The model, which we denote $\mathbf{D.6}_{-1}^{3+}$, is
		\begin{gather*}
			\boxed{
 \mfh^{\phi} \cong \mfso_{1, 3} , \qquad
				\mfk^{\phi} = \Span{\mbC} , \qquad
				\mfd^{\phi} = \Span{\mbA + \mbD_{\mbA}, \mbB + \mbD_{\mbB}} \oplus \mfk^{\phi} .
			}
		\end{gather*}
\end{itemize}

\subsection[The real forms of model $\mathbf{D.6}_{\infty}$]{The real forms of model $\boldsymbol{\mathbf{D.6}_{\infty}}$}\label{subsection:D.6-infty-real}

The adapted basis is
\begin{alignat*}{5}
& e_1 := \mbX, \qquad && e_2 := \mbY, \qquad && e_3 := \mbH, \qquad && e_4 := \mbX + \mbV_1, & \\
& e_5 := \mbY + \mbV_2, \qquad && e_6 := \mbH + 2 \mbZ . \qquad && && &
\end{alignat*}
Proceeding as for models $\mathbf{D.6}_{\lambda}$ in Section~\ref{subsection:D.6-lambda-real} gives that $\phi$ has the form \eqref{equation:phi-D.6-lambda}. Then, $\sigma_{16} = \sigma_{36} = 0$ implies that $\gamma = \zeta = \pm 1$, and then $\sigma_{34} = \sigma_{35} = 0$ gives $\alpha_1 = \pm \delta_4$, $\alpha_2 = \mp \delta_5$, $\beta_1 = \mp \epsilon_4$, $\beta_2 = \pm \epsilon_5$. Conjugating any such automorphism by another fixes the sign $\pm$, so anti-involutions with differing signs $\pm$ cannot be equivalent.

\smallskip
\noindent \textbf{Case: $\zeta = + 1$.} In this case, $\sigma_{46} = \sigma_{56} = 0$ implies that $\delta_5 = \epsilon_4 = 0$, and then $\sigma_{15} = 0$ implies that $\epsilon_5 = \delta_4^{-1}$, so $\phi(e_1) = \delta e_1$, $\phi(e_2) = \delta^{-1} e_2$, $\phi(e_3) = e_3$, $\phi(e_4) = \delta e_4$, $\phi(e_5) = \delta^{-1} e_5$, $\phi(e_6) = e_6$, where $\delta := \delta_4$, and the anti-involution condition is $|\delta| = 1$. The admissible automorphism $(e_1, \ldots, e_6) \mapsto \big(t e_1, t^{-1} e_2, e_3, t e_4, t^{-1} e_5, e_6\big)$ induces $\delta \rightsquigarrow t^2 \delta / |t|^2$, so we may normalize to $\delta = 1$, for which $\mfh^{\phi} = \spanOperator_{\bbR}\{e_1, \ldots, e_6\}$. The model, which we denote $\mathbf{D.6}_{\infty}^1$, is
\begin{gather*}
	\boxed{\begin{split}
		&\mfh^{\phi} = \mfsl_2 \oplus \big(\mfso_{1, 1} \rightthreetimes \bbR^{1, 1}\big) , \qquad
		\mfk^{\phi} = \Span{\mbH + 2 \mbZ} , \\
	&	\mfd^{\phi} = \Span{\mbX + \mbV_1, \mbY + \mbV_2} \oplus \mfk^{\phi} .
	\end{split}}
\end{gather*}

\noindent \textbf{Case: $\zeta = - 1$.} Proceeding as in the case $\zeta = + 1$ gives $\delta_4 = \epsilon_5 = 0$, $\epsilon_4 = \delta_5^{-1}$, so that $\phi(e_1) = \delta e_2$, $\phi(e_2) = \delta^{-1} e_1$, $\phi(e_3) = -e_3$, $\phi(e_4) = \delta e_5$, $\phi(e_5) = \delta^{-1} e_4$, $\phi(e_6) = -e_6$, where $\delta := \delta_5$, and the anti-involution condition is $\Im \delta = 0$. The admissible automorphism $(e_1, \ldots, e_6) \mapsto \big(t e_2, t^{-1} e_1, -e_3, t e_5, t^{-1} e_4, -e_6\big)$ induces $\delta \rightsquigarrow |t|^2 / \delta$, so we may normalize to $\delta = \pm 1$. The two choices of sign determine fixed point Lie algebras whose Killing forms have different signatures, so the resulting algebraic models are inequivalent.
\begin{itemize}[leftmargin=\caseIndent]
	\item[] \textbf{Subcase: $\delta = 1$.} In this case $\mfh^{\phi} = \spanOperator_{\bbR}\{e_1 + e_2, {\rm i} (e_1 - e_2), {\rm i} e_3, e_4 + e_5, {\rm i} (e_4 - e_5), {\rm i} e_6\}$. We may identify $\mfh^{\phi} \cong \mfsl_2 \oplus \big(\mfso_2 \rightthreetimes \bbR^2\big)$ via
		\begin{alignat*}{3}
			 &\mbX\leftrightarrow \tfrac{1}{2 \sqrt{2}} [(1 + {\rm i}) e_1 + (1 - {\rm i}) e_2] + \tfrac{1}{2} {\rm i} e_3 , \qquad
			&&\hat\mbZ\leftrightarrow \tfrac{1}{2} {\rm i} (-e_3 + e_6) ,& \\
			& \mbY\leftrightarrow \tfrac{1}{2 \sqrt{2}} [(1 + {\rm i}) e_1 + (1 - {\rm i}) e_2] - \tfrac{1}{2} {\rm i} e_3 , \qquad
			&& \mbV_1\leftrightarrow e_1 + e_2 - e_4 - e_5 ,& \\
			& \mbH\leftrightarrow \tfrac{1}{\sqrt{2}} [(-1 + {\rm i}) e_1 + (-1 - {\rm i}) e_2] , \qquad
			&& \mbV_2\leftrightarrow {\rm i}(e_1 - e_2 - e_4 + e_5) .&
		\end{alignat*}
		Here we realize $\mfso_2$ as the Lie algebra preserving the standard inner product $\left(\begin{smallmatrix}1&\cdot\\ \cdot&1\end{smallmatrix}\right)$ on $\bbR^2$, written with respect to the basis $(\mbV_1, \mbV_2)$, and we take $\hat\mbZ$ to be its standard generator, so that its action is given by
		\begin{gather*}
			[\hat\mbZ, \mbV_1] = \mbV_2 , \qquad
			[\hat\mbZ, \mbV_2] = -\mbV_1 .
		\end{gather*}
		The model, which we denote by $\mathbf{D.6}_{\infty}^2$, is
		\begin{equation*}
			\boxed{
 \begin{split}
 &\mfh^{\phi} = \mfsl_2 \oplus \big(\mfso_2 \rightthreetimes \bbR^2\big) , \qquad
	 	 \mfk^{\phi} = \operatorname{span}\big\{\mbX - \mbY + 2 \hat\mbZ\big\} , \\
 & \mfd^{\phi} =
 \operatorname{span}\big\{
 \sqrt{2} \mbX + \sqrt{2} \mbY - \mbV_1 - \mbV_2 ,
							\sqrt{2} \mbH + \mbV_1 - \mbV_2
						\big\} \oplus \mfk^{\phi} .
 \end{split}
			}
		\end{equation*}
	\item[] \textbf{Subcase: $\delta = -1$.} In this case $\mfh^{\phi} = \spanOperator_{\bbR}\{e_1 - e_2, {\rm i} (e_1 + e_2), {\rm i} e_3, e_4 - e_5, {\rm i} (e_4 + e_5), {\rm i} e_6 \}$. We may identify $\mfh^{\phi} \cong \mfso_3 \oplus \big(\mfso_2 \rightthreetimes \bbR^2\big)$ via
		\begin{alignat*}{3}
			 & \mbA\leftrightarrow \tfrac{1}{2} ( e_1 - e_2) , \qquad
			&&\hat\mbZ\leftrightarrow \tfrac{1}{2} {\rm i} (- e_3 + e_6) ,& \\
			 & \mbB\leftrightarrow \tfrac{1}{2} {\rm i} ( e_1 + e_2) , \qquad
			&& \mbV_1\leftrightarrow \tfrac{1}{2} (- e_1 + e_2 + e_4 - e_5) , &\\
			& \mbC\leftrightarrow \tfrac{1}{2} {\rm i} e_3 , \qquad
			&& \mbV_2\leftrightarrow \tfrac{1}{2} {\rm i} (- e_1 - e_2 + e_4 + e_5) .&
		\end{alignat*}
		The model, which we denote by $\mathbf{D.6}_{\infty}^4$, is
		\begin{gather*}
			\boxed{\begin{split}
& \mfh^{\phi} = \mfso_3 \oplus \big(\mfso_2 \rightthreetimes \bbR^2\big), \qquad
				\mfk^{\phi} = \operatorname{span}\big\{\mbC + \hat\mbZ\big\}, \\
& \mfd^{\phi} = \Span{
								\mbA \!+\! \mbV_1 ,
								\mbB \!+\! \mbV_2
							 } \oplus \mfk^{\phi} .
\end{split}}
		\end{gather*}
\end{itemize}

\subsection{The real forms of model \texorpdfstring{$\mathbf{D.6}_{\ast}$}{D.6-star}}\label{subsection:D.6-ast-real}
The adapted basis is
\begin{alignat*}{3}
	 e_1 &:= -\mbW_1 + \mbW_2 , \qquad
	&e_2 &:= \mbW_3 , \qquad
	&e_3 &:= \mbW_1 + \mbW_2 , \\
	 e_4 &:= \mbH - 2 \mbW_3 , \qquad
	&e_5 &:= \mbX - \mbY + 2 \mbW_1 - 2 \mbW_2 , \qquad
	&e_6 &:= \mbX + \mbY .
\end{alignat*}
We have $\mfe := [\mfk, \mfd] = \Span{e_4, e_5}$, and if we denote $\mff := \mfrad(\mfh) \cap (\mfd + [\mfd, \mfd]) = \Span{e_3}$, then $[\mfe, \mff] = \Span{e_1, e_2}$, so $\phi$ has the form \eqref{equation:phi-D.6-lambda}.

Now, $\sigma_{34} = \sigma_{35} = 0$ implies that $\alpha_1 = \gamma \delta_4$, $\alpha_2 = \gamma \delta_5$, $\beta_1 = \gamma \epsilon_4$, $\beta_2 = \gamma \epsilon_5$, then $\sigma_{16} = \sigma_{26} = 0$ gives $\zeta = \pm 1$, $\epsilon_4 = \pm \delta_5$, and $\epsilon_5 = \pm \delta_4$. Next, $\sigma_{45} = 0$ implies that $\gamma = \zeta$ and $\delta_4^2 - \delta_5^2 = 1$, so we may parameterize $\delta_4 = \cosh \tau$, $\delta_5 = \sinh \tau$, and every map satisfying the given conditions is an automorphism. Conjugating any such automorphism of this form by another fixes the sign $\pm$, so admissible anti-involutions with differing signs $\pm$ cannot be equivalent. In summary, we have
\begin{alignat*}{3}
	 \phi(e_1) &= \pm (\cosh \tau \, e_1 + \sinh \tau \, e_2) , \qquad
	&\phi(e_2) &= \sinh \tau \, e_1 + \cosh \tau \, e_2 , \qquad
	&\phi(e_3) &= \pm e_3 , \\
	 \phi(e_4) &= \cosh \tau \, e_4 + \sinh \tau \, e_5 , \qquad
	&\phi(e_5) &= \pm (\sinh \tau \, e_4 + \cosh \tau \, e_5) , \qquad
	&\phi(e_6) &= \pm e_6 .
\end{alignat*}

\noindent \textbf{Case: $\zeta = +1$.} The anti-involution condition is $\Re\tau = 0$. The admissible automorphism $(e_1, \ldots, e_6) \!\mapsto\! (\cosh t \, e_1 + \sinh t \, e_2, \sinh t \,e_1 + \cosh t \,e_2, e_3, \cosh t \, e_4 + \sinh t \, e_5, \sinh t \, e_4 + \cosh t \, e_5, e_6)$ induces $\tau \mapsto \tau + 2 {\rm i} \Im t$, so we may normalize to $\tau = 0$, for which $\mfh^{\phi} = \spanOperator_{\bbR}\{e_1, \ldots, e_6\}$. The corresponding model, which we denote $\mathbf{D.6}_{\ast}^{1-}$, is
\begin{equation*}
 \boxed{
 \begin{split}
& \mfh^{\phi} = \mfsl_2 \rightthreetimes \bbR^{1, 2} , \qquad
	 	 \mfk^{\phi} = \Span{\mbX + \mbY} , \\
& \mfd^{\phi} = \Span{\mbX - \mbY + 2 \mbW_1 - 2 \mbW_2, \mbH - 2 \mbW_3} \oplus \mfk^{\phi} ;
 \end{split}
	}
\end{equation*}
$\mfh^{\phi}$ is the affine Lorentzian algebra.

\smallskip
\noindent \textbf{Case: $\zeta = -1$.} The anti-involution condition is that $\Im \tau = k \pi$ for some integer $k$. The admissible automorphism $(e_1, \ldots, e_6) \rightsquigarrow (-\cosh t \, e_1 + \sinh t \, e_2, -\sinh t \, e_1 + \cosh t \, e_2, -e_3, \cosh t \, e_4 - \sinh t \, e_5, \sinh t \, e_4 - \cosh t \, e_5, -e_6 )$ induces $\tau \rightsquigarrow \tau + 2 \Re t$, so we may normalize to $\Re \tau = 0$.

\begin{itemize}[leftmargin=\caseIndent]
	\item[] \textbf{Subcase: $k$ even.} In this case, $\mfh^{\phi} = \spanOperator_{\bbR}\{{\rm i} e_1, e_2, {\rm i} e_3, e_4, {\rm i} e_5, {\rm i} e_6\}$. We may identify $\mfh^{\phi} \cong \mfsl_2 \rightthreetimes \bbR^3$ (where the semidirect product is as in the case $\zeta = + 1$) via
	\begin{alignat*}{4}
		 \mbX &\leftrightarrow {\rm i}\big(e_1 + \tfrac{1}{2} e_5 + \tfrac{1}{2} e_6\big) , \qquad
		&\mbW_1 &\leftrightarrow {\rm i}\big(\tfrac{1}{2} e_1 - \tfrac{1}{2} e_3\big) , \\
		 \mbY &\leftrightarrow {\rm i}\big(e_1 + \tfrac{1}{2} e_5 - \tfrac{1}{2} e_6\big) , \qquad
		&\mbW_2 &\leftrightarrow {\rm i}\big(\tfrac{1}{2} e_1 + \tfrac{1}{2} e_3\big) , \\
		 \mbH &\leftrightarrow 2 e_2 + e_4 , \qquad
		&\mbW_3 &\leftrightarrow -e_2 .
	\end{alignat*}
The model, which we denote $\mathbf{D.6}_{\ast}^{1+}$, is
		\begin{equation*}
			\boxed{
 \begin{split}
 &\mfh^{\phi} = \mfsl_2 \rightthreetimes \bbR^3 , \qquad
	 	 \mfk^{\phi} = \Span{\mbX - \mbY} , \\
 &\mfd^{\phi} = \Span{\mbH + 2 \mbW_3, \mbX + \mbY - 2 \mbW_1 - 2 \mbW_2} \oplus \mfk^{\phi} .
 \end{split}
			}
		\end{equation*}
	\item[] \textbf{Subcase: $k$ odd.} In this case $\mfh^{\phi} = \spanOperator_{\bbR}\{e_1, {\rm i} e_2, {\rm i} e_3, {\rm i} e_4, e_5, {\rm i} e_6\}$. We may identify the real form with the real Euclidean algebra, $\mfh^{\phi} = \mfso_3 \rightthreetimes \bbR^3$. Take the standard basis $(\mbW_{\mbA}, \mbW_{\mbB}, \mbW_{\mbC})$ of $\bbR^3$, so that the action is characterized by $[\mbA, \mbW_{\mbB}] = \mbW_{\mbC}$, its cyclic permutations in $\mbA, \mbB, \mbC$, and $[\mbA, \mbW_{\mbA}] = [\mbB, \mbW_{\mbB}] = [\mbC, \mbW_{\mbC}] = 0$. Then, we may identify
	\begin{alignat*}{4}
		 \mbA &\leftrightarrow e_1 + \tfrac{1}{2} e_5 , \qquad
		&\mbW_{\mbA} &\leftrightarrow e_1 , \\
		 \mbB &\leftrightarrow {\rm i} \big(e_2 + \tfrac{1}{2} e_4\big) , \qquad
		&\mbW_{\mbB} &\leftrightarrow {\rm i} e_2 , \\
		 \mbC &\leftrightarrow -\tfrac{1}{2} {\rm i} e_6 , \qquad
		&\mbW_{\mbC} &\leftrightarrow {\rm i} e_3 .
	\end{alignat*}
The model, which we denote $\mathbf{D.6}_{\ast}^3$, is
	\begin{gather*}
		\boxed{
 \mfh^{\phi} = \mfso_3 \rightthreetimes \bbR^3 , \qquad
			\mfk^{\phi} = \Span{\mbC} , \qquad
			\mfd^{\phi} = \Span{\mbA - \mbW_{\mbA}, \mbB - \mbW_{\mbB}} \oplus \mfk^{\phi} .
		}
	\end{gather*}
\end{itemize}

This completes the real classification.

\begin{Example}\label{example:D.6-ast-coordinates}
We indicate, using $\mathbf{D.6}_{\ast}^3$ as an example, how to realize a multiply transitive homogeneous distribution in local coordinates starting from an algebraic model.

First, choose a group realizing $\mfh^{\phi}$, say,
\begin{gather*}
	H
		:= \SO(3, \bbR) \rightthreetimes \bbR^3
		=
	\left\{
		\begin{pmatrix}
			1 & 0 \\
			w & R
		\end{pmatrix}
	\colon
		R \in \SO(3, \bbR),\, w \in \bbR^3
	\right\} .
\end{gather*}
We parameterize $\SO(3, \bbR)$ (and hence $H$) explicitly using Euler angles: Let $R_{\mbA}(\lambda) \in \SO(3, \bbR)$ denote the anticlockwise rotation about the oriented $\mbW_{\mbA}$-axis through an angle $\lambda$, and define $R_{\mbB}(\mu), R_{\mbC}(\nu)$ analogously. Then, for appropriate restrictions on the domain, $R(\lambda, \mu, \nu) = R_{\mbA}(\lambda) R_{\mbB}(\mu) R_{\mbC}(\nu)$ defines local coordinates $(\lambda, \mu, \nu)$ on $\SO(3, \bbR)$ and, using standard coordinates $(r, s, t)$ on $\bbR^3$, coordinates $(\lambda, \mu, \nu, r, s, t)$ on~$H$. Reading a left-invariant coframe from appropriate components of the Maurer-Cartan form
\begin{gather*}
	\begin{pmatrix}
		1 & \cdot \\
		w & R
	\end{pmatrix}^{-1}
	d
	\begin{pmatrix}
		1 & \cdot \\
		w & R
	\end{pmatrix}
 \in \Gamma(T^*H \otimes \mfh),
 \qquad w = (r, s, t)^{\top} ,
\end{gather*}
and forming the dual, left-invariant frame gives the identifications
\begin{gather*}
	\mbA\leftrightarrow
		\sec \mu \cos \nu \,\partial_{\lambda}
			+ \sin \nu \, \partial_{\mu}
			- \tan \mu \cos \nu \,\partial_{\nu} , \\
	\mbB\leftrightarrow
		-\sec \mu \sin \nu \,\partial_{\lambda}
			+ \cos \nu \,\partial_{\mu}
			+ \tan \mu \sin \nu \,\partial_{\nu} , \\
	\mbC\leftrightarrow
		\partial_{\nu} , \\
	\mbW_{\mbA}\leftrightarrow
		\cos \mu \cos \nu \, \partial_r
			+ (\cos \lambda \sin \nu + \sin \lambda \sin \mu \cos \nu) \partial_s \\
\hphantom{\mbW_{\mbA}\leftrightarrow}{} + (\sin \lambda \sin \nu - \cos \lambda \sin \mu \cos \nu) \partial_t , \\
	\mbW_{\mbB}\leftrightarrow
		-\cos \mu \sin \nu \partial_r
			+ (\cos \lambda \cos \nu - \sin \lambda \sin \mu \sin \nu) \partial_s \\
\hphantom{\mbW_{\mbB}\leftrightarrow}{} + (\sin \lambda \cos \nu + \cos \lambda \sin \mu \sin \nu) \partial_t , \\
	\mbW_{\mbC}\leftrightarrow
		\sin \mu \,\partial_r
			- \sin \lambda \cos \mu \,\partial_s
			+ \cos \lambda \cos \mu \,\partial_t .
\end{gather*}
In these coordinates the fibers of $H \to H / K$ are the integral curves of $\partial_{\nu}$, so we may use $(\lambda, \mu, r, s, t)$ as coordinates on the quotient space. Pulling back the $1$-forms defining $\mfd$ (viewed under this identification as a local distribution on $H$) by a suitable local section and computing the annihilator gives a coordinate expression for the distribution $\mbD$:
\begin{gather*}
	\operatorname{span}\big\{
		\partial_{\lambda}
			- \cos^2 \mu \, \partial_r
			- \sin \lambda \sin \mu \cos \mu \, \partial_s
			+ \cos \lambda \sin \mu \cos \mu \, \partial_t ,
		\partial_{\mu}
			- \cos \lambda \partial_s
			- \sin \lambda \partial_t
	\big\} .
\end{gather*}
Following the procedure in \cite[Section~2]{Strazzullo} allows us to put this distribution in preferred forms convenient for other purposes. In coordinates $(x^i)$,
\begin{gather*}
	x^1 = \lambda, \qquad
	x^2 = s \cos \lambda + t \sin \lambda + \mu, \qquad
	x^3 = -s \sin \lambda + t \cos \lambda, \\
	x^4 = r, \qquad
	x^5 = \tan \mu ,
\end{gather*}
$\mbD$ is the common kernel of the $1$-forms
\begin{gather*}
	-x^3 \,{\rm d}x^1 + {\rm d}x^2, \qquad
	- f \,{\rm d}x^1 - x^5 \,{\rm d}x^3 + {\rm d}x^4, \qquad
	- (\partial_{x^5} f) \,{\rm d}x^1 + {\rm d}x^3 ,
\end{gather*}
where $f := x^5 \arctan x^5 - x^2 x^5 + 1$; distributions in this form for some function $f\big(x^1, \ldots, x^5\big)$ are said to be in \textit{Goursat normal form}. For any distribution in that form, changing to coordinates $(x, y, p, q, z)$, where
\begin{gather*}
	x = x^1, \qquad
	y = x^2, \qquad
	p = x^3, \qquad
	q = \partial_{x^5} f, \qquad
	z = x^4 ,
\end{gather*}
realizes \looseness=1 the distribution in Monge normal form with the function $F$ given by writing $x^5 \partial_{x^5} f - f$ in the coordinates $(x, y, p, q, z)$. In our case, $\partial_{x^5} f = x^5 / \big[\big(x^5\big)^2 + 1\big] + \arctan x^5 - x^2$, and so one Monge normal form for $\mbD$ is given by $F(x, y, p, q, z) = -\cos^2 g^{-1}(y + q)$; here $g$ is the map $u \mapsto u + \sin u \cos u$, which, up to composition with appropriate affine transformations, appears in Kepler's equation and in a standard parametrization of the cycloid (cf.\ \cite[Section~4]{CartanFiveVariables}).
\end{Example}

\section{Identification algorithms}\label{section:identification}
We now present an algorithm for identifying the isomorphism type of a given algebraic model $(\mfh, \mfk; \mfd)$ with an explicit model in the classification; this amounts to generating sufficiently many invariants to distinguish different distributions. We split cases according to $\dim \mfh$; up to isomorphism there is only one model with $\dim \mfh = 14$ in both the complex and real settings, leaving the cases $\dim \mfh = 7, 6$.

To validate these algorithms it suffices to verify that they identify correctly each algebraic model in the corresponding (complex or real) classifications.

\subsection[Models with $\dim \mfh = 7$]{Models with $\boldsymbol{\dim \mfh = 7}$}

In this case, in both the complex and real settings, the models are determined up to equivalence by the underlying Lie algebras $\mfh$, so it is enough to distinguish those.

\subsubsection{Complex models}\label{subsubsection:identification-submaximal-complex}
(Motivated by the discussion in \cite[Section~2]{DoubrovKruglikov},) consider the maps $\ad \mbv$, $\mbv \in \mfh$. Not all maps $\ad \mbv$ are tracefree, so $\mft := \{ \mbv \in \mfh \colon \tr \ad \mbv = 0 \}$ has dimension $6$, but the nilradical $\mfN := \mfnil(\mfh) < \mft$ of $\mfh$ is isomorphic to $\mfn \otimes \bbC$ and so has dimension~$5$. Thus, the spectrum of $\ad \mbv_{\circ}$, $\mbv_{\circ} \in \mft$, depends only on the projection of $\mbv_{\circ}$ to the line $\mft / \mfN$, and so the conformal spectrum of $\ad \mbv_{\circ}$ is independent of the choice of $\mbv_{\circ} \in \mft - \mfN$. In particular, where $\sigma_k$ denotes the $k$th symmetric polynomial in the eigenvalues of $\ad \mbv_{\circ}$,
\begin{gather*}
	\boxed{\Lambda := \frac{64 \sigma_4}{100 \sigma_4 - 9 \sigma_2^2}}
\end{gather*}
is an invariant of the conformal spectrum of $\ad \mbv_{\circ}$ and hence an invariant of $\mfh$.\footnote{The invariant $J$ for these distributions defined in \cite{DoubrovKruglikov} is $4 \sigma_4 / \sigma_2^2 = 9 \Lambda / (25 \Lambda - 16)$; we have $\Lambda = 16 J / (25 J - 9)$.}

For the Lie algebras $\mfh_{r, s}^{\bbC} = (\mfn \otimes \bbC) \rightthreetimes \Span{\mbE, \mbF}$ of the submaximal complex algebraic models, $\mfN = \mfn \otimes \bbC$ and $\mft = (\mfn \otimes \bbC) \oplus \Span{\mbF}$, so we need only compute the spectrum of $\ad \mbF$. Consulting the formulae in Section~\ref{subsection:N.7} gives that the spectrum is $(+a, -a, +b, -b, 0, 0, 0)$, so $\sigma_2 = -\big(a^2 + b^2\big) = -r$ and $\sigma_4 = a^2 b^2 = s$, and thus $\Lambda = 64 s / \big(100 s - 9 r^2\big)$. Specializing to Cartan's Monge normal form ($r = \frac{10}{3} I$, $s = I^2 + 1$) gives $\Lambda = 16 \big(I^2 + 1\big) / 25$; thus, since $I^2$ is a~complete invariant, so is $\Lambda$.

\subsubsection{Real models}\label{subsubsection:identification-submaximal-real}

It follows from the form of the equivalence relation in Section~\ref{subsection:N.7-real} that every submaximal complex model that admits a real form admits precisely two. So, to identify a submaximal real algebraic model $(\mfh, \mfk; \mfd)$, we first identify its complexification $(\mfh \otimes \bbC, \mfk \otimes \bbC; \mfd \otimes \bbC)$.

If we define $\mft, \mfN$ as in the complex case, then replacing a choice $\mbv_{\circ} \in \mft - \mfN$ with $t \mbv_{\circ} + \mathbf{n}$, $t \in \bbR^*$, $\mathbf{n} \in \mfN$ induces $\sigma_k \rightsquigarrow t^k \sigma_k$, so the signs of $\sigma_2$, $\sigma_4$ are independent of the choice $\mbv_{\circ}$. For real pairs $(r, s)$, the formulae derived in the complex case give for $\mfh_{r, s}$ and the choice $\mbv_{\circ} = \mbF$ that $\sign(\sigma_2) = -\sign(r)$ and $\sign(\sigma_4) = \sign(s)$, so for all $(r, s)$ the pair $(\sign(\sigma_2), \sign(\sigma_4))$ distinguishes between the real forms.

\subsection[Models with $\dim \mfh = 6$]{Models with $\boldsymbol{\dim \mfh = 6}$}\label{subsection:identification-6}

\subsubsection{Complex models}\label{subsubsection:identification-6-complex} The isomorphism type of $\mfh$ can be determined by examining its radical, $\mfr := \mfrad(\mfh)$: Consulting Table~\ref{table:complex-6} shows that if $\mfr \neq \{ 0 \}$, then $\dim \mfr = 3$, and the derived algebra $[\mfr, \mfr]$ has dimension $0$, $1$, or $2$, respectively, if $\mfh$ is isomorphic to $\mfso_3(\bbC) \rightthreetimes \bbC^3$, $\mfsl_2(\bbC) \rightthreetimes (\mfm \otimes \bbC)$, or $\mfsl_2(\bbC) \oplus \big(\mfso_2(\bbC) \rightthreetimes \bbC^2\big)$, in which case the model is isomorphic, respectively, to $\mathbf{D.6}_{\ast}$, $\mathbf{N.6}$, or $\mathbf{D.6}_{\infty}$.

If instead $\mfr = \{ 0 \}$, then $\mfh \cong \mfsl_2(\bbC) \oplus \mfsl_2(\bbC)$, and the model is isomorphic to $\mathbf{D.6}_{\lambda}$ for some parameter $\lambda$. Let $\pi_1, \pi_2 \colon \mfh \to \mfsl_2(\bbC)$ denote the projections onto the two summands, and define $\mfe := [\mfk, \mfd]$. Then, if $Q$ is the quadratic form on $\mfsl_2(\bbC)$ induced by the Killing form, the restriction of $\pi_2^* Q$ to $\mfe$ is some constant multiple of the restriction of $\pi_1^* Q$ to~$\mfe$. Reversing the assignment of the indices $1$, $2$ to the summands $\mfsl_2(\bbC)$ replaces the constant with its reciprocal, so up to this inversion the constant is an invariant of the algebraic model.

For $\mathbf{D.6}_{\lambda}$, $\mfe = \Span{\mbX - \lambda \mbX', \mbY - \mbY'}$, so for $\mbv := u(\mbX -\lambda \mbX') + v (\mbY - \mbY') \in \mfe$, $(\pi_1^* Q)(\mbv) = Q(u \mbX + v \mbY) = 8 u v$ and $(\pi_2^* Q)(\mbv) = Q(-u \lambda \mbX' - v \mbY') = 8 \lambda u v$, and so the constant is $\lambda$.

Applying this algorithm to Cartan's realization of these distributions in \cite[Section~50]{CartanFiveVariables} gives that the invariant $\lambda$ can be written (again, up to inversion) in terms of his inva\-riants~$\rho$ and~$\sigma$ (defined in Section~49) as $(3 \sigma + 7 \rho) / (3 \rho + 7 \sigma)$ (cf. the third display equation in Section~50 of that reference).

\subsubsection{Real models}\label{subsubsection:identification-6-real}
As in the case $\dim \mfh = 7$, we first classify the complexification of the model, which reduces the real classification to distinguishing among the finitely many real forms thereof.

\begin{itemize}[leftmargin=1.5cm]\itemsep=0pt
\item[\textbf{$\mathbf{N.6}$}] We construct a basis $(f_a)$ of $\mfh$ adapted to the model: Fix $f_6 \in \mfk - \{ 0 \}$, pick an element $f_5 \in [\mfk, \mfd]$ so that $[f_5, f_6] = 2 f_6$, and choose $f_4$ so that $(f_4, f_5, f_6)$ is a basis of $\mfd$.

As a vector space (not as a Lie algebra) $\mfh$ decomposes as a direct sum $\mfrad(\mfh) \oplus \mfd$; denote by $f_3$ the projection of $[f_4, f_5]$ onto the first summand. Then, $f_1 := [f_3, f_4]$ and $f_2 := [f_3, f_5] + f_3$ satisfy $[f_1, f_3] = \mu f_2$ for some $\mu \neq 0$; $\mu$ is independent of the choices of $f_5$ and $f_4$, and replacing $f_6$ with $t f_6$, $t \neq 0$, induces $\mu \rightsquigarrow t^2 \mu$, so $\sign(\mu) = \pm 1$ is an invariant of the model, and the model is respectively isomorphic to $\mathbf{N.6}^{\pm}$.

\item[\textbf{$\mathbf{D.6}_{\lambda}$}] The isomorphism type of $\mfh$ can be distinguished by its Killing form $\kappa$ (see Table \ref{table:real-forms-sl2C-sl2C}). If $\mfh \cong \mfso_3 \oplus \mfsl_2$ or $\mfh \cong \mfso_3 \oplus \mfso_3$, then the model is isomorphic to $\mathbf{D.6}_{\lambda}^4$ or $\mathbf{D.6}_{\lambda}^6$, respectively. In both remaining cases there are two real forms, and they can be distinguished by the restriction of $\kappa$ to $\mfk$: If $\mfh \cong \mfsl_2 \oplus \mfsl_2$, then the model is isomorphic to $\mathbf{D.6}_{\lambda}^{2-}$ or $\mathbf{D.6}_{\lambda}^{2+}$ if $\kappa\vert_{\mfk}$ is positive or negative definite, respectively. If instead $\mfh \cong \mfso_{1, 3}$ (which occurs only for $\lambda = -1$), then the model is isomorphic to $\mathbf{D.6}_{-1}^-$ or $\mathbf{D.6}_{-1}^+$ if $\kappa\vert_{\mfk}$ is positive or negative definite, respectively.

\item[\textbf{$\mathbf{D.6}_{\infty}$}] The three possibilities are distinguished by the isomorphism type of $\mfh$, or just as well by the signature of the Killing form of $\mfh$: For each model $\mathbf{D.6}_{\infty}^s$ the signature is $(4 - s, s)$.

\item[\textbf{$\mathbf{D.6}_{\ast}$}] If the Killing form of $\mfh$ is definite the model is isomorphic to $\mathbf{D.6}_{\ast}^3$. Otherwise, $\mfh \cong \mfsl_2 \rightthreetimes \bbR^3$, in which case the two possibilities are again distinguishable via restriction of the Killing form $\kappa$ to $\mfk$: The model is isomorphic to $\mathbf{D.6}_{\infty}^{1-}$ or $\mathbf{D.6}_{\infty}^{1+}$ if $\kappa\vert_{\mfk}$ is positive- or negative-definite, respectively.
\end{itemize}

\begin{Example}\label{example:N.6-Monge}
Consider the real distribution defined by $F(x, y, p, q, z) := q^{1 / 3} - y$. Replacing the line \texttt{F := q\^{}(1 / 3) + y;} in the program in Example \ref{example:coordinates-to-distribution} with \texttt{F := q\^{}(1 / 3) - y;} and executing gives that $\mfh := \mfaut(\mbD_F)$ has basis $(\xi_i')$, where $\xi_1' := -y \partial_x + p^2 \partial_p + 3 p q \partial_q + \frac{1}{2} y^2 \partial_z$, $\xi_2' = -\big(x \partial_y + \partial_p - \frac{1}{2} x^2 \partial_z\big)$, $\xi_3' = \xi_3$, $\xi_4' = \xi_4$, $\xi_5' = \xi_5$, $\xi_6' = -\partial_y + x \partial_z$, and where the $\xi_i$ were defined in \eqref{equation:N.6-symmetries}; in particular, the distribution is locally homogeneous. We have $\mfr := \mfrad(\mfh) = \Span{\xi_4', \xi_5', \xi_6'}$, and $[\mfr, \mfr] = \Span{\xi_4'}$ has dimension~$1$, so the distribution is (locally equivalent to) a real form of $\mathbf{N.6}$.

At the base point $(0, 0, 0, 1, 0)$, the isotropy subalgebra is $\mfk = \Span{\xi_1'}$, so we may take $f_6 := \xi_1'$, and $\mfd = \Span{\xi_1', \xi_2' - \xi_4' - \xi_5', \xi_3'}$. Since $\xi_3' \in [\mfk, \mfd]$ and $[\xi_3', f_6] = 2 f_6$ we may take $f_5 := \xi_3'$ and $f_4 = 2(\xi_2' - \xi_4' - \xi_5')$, and this determines $f_3 = 4 \xi_4' + 6 \xi_5'$, $f_1 = 12 \xi_6'$, and $f_2 = 4 \xi_4'$. Then, $[f_1, f_3] = -72 \xi_4' = -18 f_2$, and this last coefficient is negative, so the model is isomorphic to~$\mathbf{N.6}^-$.
\end{Example}

\section{Realizations as rolling distributions}\label{section:rolling}

An important subclass of $(2, 3, 5)$ distributions are the so-called \textit{rolling distributions}, which are defined kinematically by a system of two Riemannian surfaces $(\Sigma_1, g_1)$, $(\Sigma_2, g_2)$ rolling along one another. Following~\cite{BryantHsu}, define the ($5$-dimensional) configuration space $C$ to be the set of triples $(x_1, x_2, \Phi)$ for which $x_1 \in \Sigma_1$, $x_2 \in \Sigma_2$ (so that $x_1$, $x_2$ are the points of tangency on the surfaces), and $\Phi\colon T_{x_1} \Sigma_1 \to T_{x_2} \Sigma_2$ is an isometry (encoding the relative rotation of $\Sigma_1$, $\Sigma_2$ about the contact point). The projection $(x_1, x_2, \Phi) \mapsto (x_1, x_2)$ realizes $C$ as a principal $\operatorname{O}(2)$-bundle $C \to \Sigma_1 \times \Sigma_2$. A smooth path $\gamma\colon I \to C$, $\gamma(t) = (x_1(t), x_2(t), \Phi(t))$, in this space encodes a~rolling trajectory of the two surfaces along one another.

The physical \textit{no-slip} condition is that the point of contact moves with the same relative motion on each surface, that is, that $x_2'(t) = \Phi(t) \cdot x_1'(t)$. The \textit{no-twist} condition is that for any (equivalently every) parallel orthonormal frame field $(E_a(t))$ along $x_1(t)$, the frame field \mbox{$(\Phi(t) \cdot E_a(t))$} along $x_2(t)$ must also be parallel. These conditions together impose three independent linear constraints on each tangent space $T_{(x_1, x_2, \Phi)} C$, and so the space of admissible velocities $\gamma'$ is a~$2$-plane distribution $\mbD$ on $C$. Direct computation shows that the restriction of $\mbD$ to the open set of points $(x_1, x_2, \Phi) \in C$ where the respective Gaussian curvatures $\kappa_i({x_i})$ of $(\Sigma_i, g_i)$ at~$x_i$, $i = 1, 2$, are unequal is a $(2, 3, 5)$ distribution. In this case, the $3$-plane distribution $[\mbD, \mbD]$ is the space of velocities satisfying just the no-slip condition. For any positive constant~$c$, the pair $(\Sigma_1, c g_1)$, $(\Sigma_2, c g_2)$ gives rise to the same distribution $\mbD$.

It is immediate that a pair of surfaces $\Sigma_1, \Sigma_2$ with constant (and unequal) Gaussian curvature give rise to a $(2, 3, 5)$ distribution $(C, \mbD)$ with $\dim \mfaut(\mbD) \geq 6$: Since both metrics $g_i$ have constant curvature, the respective algebras $\mfaut(g_i)$ of Killing fields both have dimension $3$ and by functoriality lift to infinitesimal symmetries of $(C, \mbD)$, the space of which thus contains an isomorphic copy of $\mfaut(g_1) \oplus \mfaut(g_2)$.

Many distributions occurring in the real classification are realizable this way.
\begin{itemize}[leftmargin=1.5cm]\itemsep=0pt
	\item[$\mbO^{\bbR}$] two $2$-spheres whose Gaussian curvatures have ratio $9 : 1$,
	\item[$\mathbf{D.6}_{\lambda}^6$] two $2$-spheres whose Gaussian curvatures have ratio $\kappa_1 / \kappa_2 = \lambda > 0$,
	\item[$\mathbf{D.6}_{\lambda}^4$] a $2$-sphere and a hyperbolic plane normalized to have Gaussian curvature ratio $\lambda < 0$,
	\item[$\mathbf{D.6}_{\lambda}^{2+}$] two hyperbolic planes whose Gaussian curvatures have ratio $\lambda > 0$,
	\item[$\mathbf{D.6}_{\infty}^4$] a $2$-sphere and the Euclidean plane\footnote{This corrects a misstatement in \cite{WillseCounterexample}.},
	\item[$\mathbf{D.6}_{\infty}^{2}$] a hyperbolic plane and the Euclidean plane.
\end{itemize}
Up to the joint scaling of pairs of surfaces by a common factor and local equivalence, this exhausts all of the pairs of constant curvature Riemannian surfaces of unequal curvature.

\looseness=-1 We can realize some of the remaining real forms with $\dim \mfh = 6$ by extending our attention to rolling distributions generated by pairs of \textit{Lorentzian} surfaces: We proceed as before, but instead take $(\Sigma_1, g_1)$, $(\Sigma_2, g_2)$ to be Lorentzian surfaces, in which case $C \to \Sigma_1 \times \Sigma_2$ is a principal $\operatorname{O}(1, 1)$-bundle.
\begin{itemize}\itemsep=0pt
 \item The flat (zero-curvature) model of Lorentzian surfaces is the Lorentzian plane $\bbR^{1, 1}$, the affine real plane equipped with the flat Lorentzian metric $\bar g_{1, 1}$.
 \item A Lorentzian surface of negative constant curvature is locally isometric to $2$-dimensional \textit{anti-de Sitter space}: In $3$-dimensional Lorentzian space $\big(\bbR^{1, 2}, \bar g_{1, 2}\big)$ define for a parameter $\alpha > 0$ $({\rm AdS}_2, g_{{\rm AdS}})$ to be the hypersurface $\big\{\mbW \in \bbR^{1, 2} \colon \bar g_{1, 2}(\mbW, \mbW) = -\alpha^{-2} \big\}$ equipped with the (Lorentzian) pullback metric $g_{{\rm AdS}}$; it has scalar curvature $-\alpha^2$.
 \item A Lorentzian surface of positive constant curvature is locally isometric to $2$-dimensional \textit{de Sitter} space, namely, ${\rm AdS}_2$ equipped with the Lorentzian metric $g_{{\rm dS}} := -g_{{\rm AdS}}$, which has scalar curvature~$\alpha^2$. The algebra $\mfaut(g_{{\rm dS}}) = \mfaut(g_{{\rm AdS}})$ of Killing fields is isomorphic to $\mfsl_2 \cong \mfso_{1, 2}$.
\end{itemize}

We realize some of the remaining models as rolling distributions.
\begin{itemize}[leftmargin=1.5cm]\itemsep=0pt
	\item[$\mathbf{D.6}_{\lambda}^{2-}$] if $\lambda > 0$, two de Sitter spaces ${\rm dS}_2$ (or, just as well, two anti-de Sitter spaces ${\rm AdS}_2$) with ratio $\kappa_1 / \kappa_2 = \lambda$ of Gaussian curvatures; if $\lambda < 0$, de Sitter space and anti-de Sitter space with ratio $\lambda$ of Gaussian curvatures \item[$\mathbf{D.6}_{\infty}^1$] de Sitter (or anti-de Sitter) space and the Lorentzian plane.
\end{itemize}

One can, of course, use the above identifications as rolling distributions to write down local coordinate realizations for the corresponding models, and this can be done efficiently using the procedure in \cite[Section~2]{AnNurowski}.

\begin{Example}[two spheres]\label{example:two-spheres}
Realize the case of two spheres with curvature ratio $\lambda$ with spheres $\Sigma_1$, $\Sigma_2$ with radii $\sqrt{\lambda}$, $1$. Writing the round metrics as $\lambda^{-1} [{\rm d}\alpha \otimes {\rm d}\alpha + (\sin \alpha) {\rm d}\beta \otimes {\rm d}\beta]$ and ${\rm d}\gamma \otimes {\rm d}\gamma + (\sin \gamma) {\rm d}\zeta \otimes {\rm d}\zeta$, taking respective orthonormal frames $\big(\sqrt{\lambda} \partial_{\alpha}, \sqrt{\lambda} \csc \alpha \,\partial_{\beta}\big)$ and $(\partial_{\gamma}, \csc \gamma \,\partial_{\zeta})$, and applying the procedure gives the coordinate realization
\begin{gather*}
		\operatorname{span}\big\{
			\sqrt\lambda \,\partial_{\alpha}
				+ \cos \varphi \,\partial_{\gamma}
				+ \csc \gamma \sin \varphi \,\partial_{\zeta}
				- \cot \gamma \sin \varphi \,\partial_{\varphi}, \\
\hphantom{\operatorname{span}\big\{}{}\sqrt\lambda \csc \alpha \,\partial_{\beta}
				- \sin \varphi \,\partial_{\gamma}
				+ \csc \gamma \cos \varphi \,\partial_{\zeta}
				+ (\sqrt\lambda \cot \alpha - \cot \gamma \cos \varphi) \partial_{\varphi}
		\big\}
\end{gather*}
of $\mbD$; here, $\varphi$ is a standard coordinate on the fibers of $C \to \Sigma_1 \times \Sigma_2$.
\end{Example}

\begin{Remark}
Any two constant curvature surfaces (of the same signature) whose curvatures have ratio $9:1$ determine a locally flat rolling distribution. Up to replacing both metrics with their negatives, the possibilities are: two spheres, two copies of the hyperbolic plane, and two copies of de Sitter space.
\end{Remark}

\appendix
\appendixpage
\section[An explicit realization of $\mfg_2$]{An explicit realization of $\boldsymbol{\mfg_2}$}\label{appendix:tables}

The complex (Table \ref{table:complex-14}) and real (Table~\ref{table:real-O}) flat models respectively have symmetry algebras isomorphic to $\mfg_2(\bbC)$ and $\mfg_2$. We can realize $\mfg_2$ as the subalgebra \cite{SagerschnigWillseOperator}
\begin{gather}\label{equation:g2}
	\left\{
		\begin{pmatrix}
			- \tr A & Z & s & W^{\top} & \cdot \\
			X & A - (\tr A) \bbI & \sqrt{2} \bbJ Z^{\top} & \frac{1}{\sqrt{2}} s \bbJ & -W \\
			r & -\sqrt{2} X^{\top} \bbJ & \cdot & -\sqrt{2} Z \bbJ & s \\
			Y^{\top} & -\frac{1}{\sqrt{2}} \bbJ & \sqrt{2} \bbJ X & (\tr A) \bbI - A^{\top} & -Z^{\top} \\
			\cdot & -Y & r & -X^{\top} & \tr A
		\end{pmatrix}
			\colon
		\begin{array}{@{}c@{}}
			X, W \in \bbR^2 ; \\
			r, s \in \bbR ; \\
			Y, Z \in \big(\bbR^2\big)^*; \\
			A \in \mfgl_2
		\end{array}
	\right\} < \mfgl_7 ,
\end{gather}
where
\begin{gather*}
	\bbI := \begin{pmatrix}1&\cdot\\ \cdot&1\end{pmatrix}
		\qquad \textrm{and} \qquad
	\bbJ := \begin{pmatrix}\cdot&-1\\1&\cdot\end{pmatrix} .
\end{gather*}
Then, write
\begin{gather*}
	Y = \begin{pmatrix}Y_1 & Y_2\end{pmatrix}, \qquad
	X = \begin{pmatrix}X_1 \\ X_2\end{pmatrix}, \qquad
	A = \begin{pmatrix}A_{11} & A_{12} \\ A_{21} & A_{22}\end{pmatrix}, \\
	Z = \begin{pmatrix}Z_1 & Z_2\end{pmatrix}, \qquad
	W = \begin{pmatrix}W_1 \\ W_2\end{pmatrix},
\end{gather*}
set $\mbY_1 \in \mfg_2$ to be the matrix with $Y_1 = 1$ and all other parameters set to $0$, and define the other~$13$ basis elements analogously.

We may realize the standard representation of $\mfg_2$ as the restriction of the standard action of~$\mfgl_7$ on~$\bbR^7$. Then, the choice of real algebraic model in Table~\ref{table:real-O} realizes $\mfk = \mfq$ as the (parabolic) subalgebra fixing the line $\big\{(\ast, 0, \ldots, 0)^{\top}\big\} \subset \bbR^7$ in that representation, and the choice of complex algebraic model in Table~\ref{table:complex-14} realizes $\mfk = \mfq \otimes \bbC$ analogously.

\newpage

\begin{landscape}

\section{Classification tables}\label{appendixB}

\subsection[Multiply transitive homogeneous complex $(2, 3, 5)$ distributions]{Multiply transitive homogeneous complex $\boldsymbol{(2, 3, 5)}$ distributions}

\begin{table}[htb]
\caption{Multiply transitive homogeneous complex $(2, 3, 5)$ distributions: $\dim \mfh = 14$.}\label{table:complex-14}\vspace{1mm}

\tiny
\centering
\begin{tabular}{cccccl}
 \toprule
 model
 & symmetry algebra
 & \multicolumn{1}{c}{algebraic model}
 & Monge $F$
 & \multicolumn{1}{c}{$\mfaut(\mbD_F)$} \\
 \midrule
 $\mbO$
 & $\mfg_2(\bbC)$
 & $\begin{array}{rl}
 \mfk\colon
 &\!\!\!\!\!\! \mbA_{11}, \mbA_{12}, \mbA_{21}, \mbA_{22}, \\
 &\!\!\!\!\!\! \mbZ_1, \mbZ_2, \mbs, \mbW_1, \mbW_2 \\
	 \mfd / \mfk\colon
 &\!\!\!\!\!\! \mbX_1, \mbX_2
 \end{array}$
 & $q^2$
 & $\begin{array}{c}
 \begin{array}{c}
 \mbY_ 1 \leftrightarrow \tfrac{1}{6} \partial_y , \qquad
	 \mbY_ 2 \leftrightarrow \tfrac{1}{3} \partial_z , \qquad
 	\mbr \leftrightarrow -\tfrac{1}{2 \sqrt{2}}(x \partial_y + \partial_p) , \qquad
 \mbX_ 1 \leftrightarrow \tfrac{1}{2} x^2 \partial_y + x \partial_p + \partial_q + 2 p \partial_z , \qquad
	 \mbX_ 2 \leftrightarrow \partial_x , \\
 \mbA_{11} \leftrightarrow -(y \partial_y + p \partial_p + q \partial_q + 2 z \partial_z) , \qquad
	 \mbA_{12} \leftrightarrow -\tfrac{1}{6} x^3 \partial_y - \tfrac{1}{2} x^2 \partial_p - x \partial_q + (-2 x p + 2 y) \partial_z , \\
 \mbA_{21} \leftrightarrow -q \partial_x + (-p q + \tfrac{1}{2} z) \partial_y - \tfrac{1}{2} q^2 \partial_p - \tfrac{1}{3} q^3 \partial_z , \qquad
	 \mbA_{22} \leftrightarrow -(x \partial_x + 2 y \partial_y + p \partial_p + z \partial_z)
 \end{array} \\
 \arraycolsep=1.4pt\begin{array}{rcl}
	 \mbZ_ 1 &\leftrightarrow& (-3 x q + 4 p) \partial_x + (-3 x p q + \tfrac{3}{2} x z + 2 p^2) \partial_y + (-\tfrac{3}{2} x q^2 + \tfrac{3}{2} z) \partial_p - q^2 \partial_q - x q^3 \partial_z \\
 	\mbZ_ 2 &\leftrightarrow& -x^2 \partial_x - 3 x y \partial_y + (-x p - 3 y) \partial_p + (x q - 4 p) \partial_q - 4 p^2 \partial_z \\
	 \mbs &\leftrightarrow& \tfrac{1}{2 \sqrt{2}} [(-6 x^2 q + 16 x p - 12 y) \partial_x + (-6 x^2 p q + 3 x^2 z + 8 x p^2) \partial_y + (-3 x^2 q^2 + 6 x z + 4 p^2) \partial_p \\
			 & & \qquad\qquad + (-4 x q^2 + 4 p q + 6 z) \partial_q + (-2 x^2 q^3 + 12 p z) \partial_z] , \\
	 \mbW_ 1 &\leftrightarrow& (-x^3 q + 4 x^2 p - 6 x y) \partial_x + (- x^3 p q + \tfrac{1}{2} x^3 z + 2 x^2 p^2 - 6 y^2) \partial_y + (-\tfrac{1}{2} x^3 q^2 + \tfrac{3}{2} x^2 z + 2 x p^2 - 6 y p) \partial_p \\
			 & & \qquad\qquad+ (-x^2 q^2 + 2 x p q + 3 x z - 4 p^2) \partial_q + (- \tfrac{1}{3} x^3 q^3 + 6 x p z - 6 y z - \tfrac{8}{3} p^3) \partial_z \\
 	\mbW_ 2 &\leftrightarrow& (6 y q - 4 p^2) \partial_x + (6 y p q - 3 y z - \tfrac{8}{3} p^3) \partial_y + (3 y q^2 - 3 p z) \partial_p + (2 p q^2 - 3 q z) \partial_q + (2 y q^3 - 3 z^2) \partial_z
 \end{array}
 \end{array}$ \\
 \bottomrule
\end{tabular}

\vspace{0.75cm}

\caption{Multiply transitive homogeneous complex $(2, 3, 5)$ distributions: $\dim \mfh = 7$.}\label{table:complex-7}\vspace{1mm}

\tiny
\centering
\begin{tabular}{lcclcc}
\toprule
 model
 & $\Lambda$
 & $\mbF$-action
 & \multicolumn{1}{c}{algebraic model}
 & Monge $F$
 & \multicolumn{1}{c}{$\mfaut(\mbD_F)$} \\
\midrule
 & $\bbC - \{0, 1\}$
 & $\begin{pmatrix}
 a & \cdot & \cdot & \cdot \\
 \cdot & b & \cdot & \cdot \\
 \cdot & \cdot & -a & \cdot \\
 \cdot & \cdot & \cdot & -b
 \end{pmatrix}$
 & $\begin{array}{rl}
 \mfk\colon
 &\!\!\!\!\!\! \mbE, \tfrac{1}{a} \mbS_1 + \tfrac{1}{b} \mbS_2 - \mbT_1 + \mbT_2 \\
	 \mfd / \mfk\colon
 &\!\!\!\!\!\! \mbF, \mbS_1 + \mbS_2 + a \mbT_1 - b \mbT_2
 \end{array}$
 & $q^2 + r p^2 + s y^2$
 & $\arraycolsep=1.4pt\begin{array}{c}
 \mbU \leftrightarrow 4 \big(a^2 - b^2\big) \partial_z , \qquad
 \mbS_1 \leftrightarrow \xi( b, a) , \qquad
 \mbS_2 \leftrightarrow \xi( a, b) , \\
 \mbT_1 \leftrightarrow -\tfrac{1}{a} \xi( b, -a) , \qquad
 \mbT_2 \leftrightarrow \tfrac{1}{b} \xi( a, -b) , \\
 \mbE \leftrightarrow y \partial_y + p \partial_p + q \partial_q + 2 z \partial_z , \qquad
 \mbF \leftrightarrow \partial_x
 \end{array}$ \\
 \cmidrule{2-6}
 $\mathbf{N.7}_\Lambda$
 & $0$
 & $\begin{pmatrix}
 a & \cdot & \cdot & \cdot \\
 \cdot & \cdot & \cdot & 1 \\
 \cdot & \cdot & -a & \cdot \\
 \cdot & \cdot & \cdot & \cdot
 \end{pmatrix}$
 & $\begin{array}{rl}
 \mfk\colon
 &\!\!\!\!\!\! \mbE, -2 a \mbS_1 + \mbT_1 + 2 a \mbT_2 \\
	 \mfd / \mfk\colon
 &\!\!\!\!\!\! \mbF, - \mbS_2 + \mbT_1 + a \mbT_2
 \end{array}$
 & $q^2 + p^2$
 & $\arraycolsep=1.4pt\begin{array}{c}
 \mbU \leftrightarrow 2 \partial_z , \qquad
 \mbS_1 \leftrightarrow \frac{1}{2} \xi(0, 1) , \qquad
 \mbS_2 \leftrightarrow \partial_y , \\
 \mbT_1 \leftrightarrow \xi(0, -1) , \qquad
 \mbT_2 \leftrightarrow x \partial_y + \partial_p + 2 y \partial_z , \\
 \mbE \leftrightarrow y \partial_y + p \partial_p + q \partial_q + 2 z \partial_z , \qquad
 \mbF \leftrightarrow \partial_x
 \end{array}$ \\
 \cmidrule{2-6}
 & $1$
 & $\begin{pmatrix}
 a & \cdot & \cdot & 1 \\
 \cdot & -a & 1 & \cdot \\
 \cdot & \cdot & -a & \cdot \\
 \cdot & \cdot & \cdot & a
 \end{pmatrix}$
 & $\begin{array}{rl}
 \mfk\colon
 &\!\!\!\!\!\! \mbE, \mbS_2 + a \mbT_1 + \mbT_2 \\
	 \mfd / \mfk\colon
 &\!\!\!\!\!\! \mbF, \mbS_1 + a \mbS_2 + 2 a \mbT_2
 \end{array}$
 & $q^2 + 2 p^2 + y^2$
 & $\arraycolsep=1.4pt\begin{array}{c}
 \mbU \leftrightarrow - 8 \partial_z , \qquad
 \mbS_1 \leftrightarrow \xi(1, 1) , \qquad
 \mbS_2 \leftrightarrow \xi(1, -1) , \\
 \mbT_1 \leftrightarrow e^{-x} [ x \partial_y - (x - 1) \partial_p + (x - 2) \partial_q - 2 (x y + (-x p + y) + 2 p) \partial_z] , \\
 \mbT_2 \leftrightarrow e^{ x} [(x - 1) \partial_y + x \partial_p + (x + 1) \partial_q + 2 (x y + ( x p - 2 y) + p) \partial_z] , \\
 \mbE \leftrightarrow y \partial_y + p \partial_p + q \partial_q + 2 z \partial_z , \qquad
 \mbF \leftrightarrow \partial_x
 \end{array}$ \\
\bottomrule
\end{tabular}
\begin{minipage}{200mm}
\vspace{0.25cm}
\tiny
Here, $\xi(u, v) := e^{v x} \big[\partial_y + v \partial_p + v^2 \partial_q + 2 v\big(u^2 y + v p\big) \partial_z\big]$ .

For each $\Lambda \in \bbC$, $(r, s)$ is a pair that satisfies $\Lambda = 64 s / \big(100 s - 9 r^2\big)$, and $a$, $b$ satisfy $r = a^2 + b^2$, $s = a^2 b^2$; see Section~\ref{subsection:N.7} for details. In each case, the symmetry algebra $\mfh_{r, s}^{\bbC}$ is an extension of the complex $5$-dimensional Heisenberg algebra $\mfn \otimes \bbC = \Span{\mbU, \mbS_1. \mbS_2, \mbT_1, \mbT_2}$~-- for which the nonzero brackets are those determined by $[\mbS_1, \mbT_1] = [\mbS_2, \mbT_2] = \mbU$~-- by two commuting derivations, $\mbE$, $\mbF$: (1) $\mbE$ is the grading derivation, so $[\mbE, \mbS_1] = -\mbS_1$, $[\mbE, \mbS_2] = -\mbS_2$, $[\mbE, \mbT_1] = -\mbT_1$,$[\mbE, \mbT_2] = -\mbT_2$, $[\mbE, \mbU] = -2 \mbU$. (2) The derivation $\mbF$ is defined by $[\mbF, \mbU] = \mbzero$ and the given action of $\mbF$ on $\mfn_{-1}(\bbC)$ with respect to the basis $(\mbS_1, \mbS_2, \mbT_1, \mbT_2)$.
\end{minipage}
\end{table}

\newpage

\vspace*{\fill}

\begin{table}[htb]
\caption{Multiply transitive homogeneous complex $(2, 3, 5)$ distributions: $\dim \mfh = 6$.}\label{table:complex-6}\vspace{1mm}

\tiny
\centering
\setlength\tabcolsep{3pt}
\begin{tabular}{lcllcl}
 \toprule
 model
 & symmetry algebra structure
 & \multicolumn{1}{c}{algebraic model}
 & \multicolumn{1}{c}{adapted basis}
 & Monge $F$
 & \multicolumn{1}{c}{$\mfaut(\mbD_F)$} \\
 \midrule
 $\mathbf{N.6}$
 & $\begin{array}{c}
 \mfsl_2(\bbC) \rightthreetimes \mfm \otimes \bbC \\
 \begin{array}{c|cccccc}
 & \mbX & \mbY & \mbH & \mbU & \mbS & \mbT \\
	 \hline
 \mbX & \cdot & \mbH & -2 \mbX & \cdot & \cdot & \mbS \\
 \mbY & & \cdot & 2 \mbY & \cdot & \mbT & \cdot \\
 \mbH & & & \cdot & \cdot & \mbS & -\mbT \\
 \mbU & & & & \cdot & \cdot & \cdot \\
 \mbS & & & & & \cdot & \mbU \\
 \mbT & & & & & & \cdot \\
 \end{array}
 \end{array}$
 & $\begin{array}{rl}
 \mfk\colon
 &\!\!\!\!\!\! \mbX \\
	 \mfd / \mfk\colon
 &\!\!\!\!\!\! \mbY - \mbU - \mbS, \\
 &\!\!\!\!\!\! \mbH
 \end{array}$
 & $\arraycolsep=1.4pt\begin{array}{rcl}
 e_1 &=& 3 \mbT \\
 e_2 &=& \mbU \\
 e_3 &=& 2 \mbU + 3 \mbS \\
 e_4 &=& \mbY - \mbU - \mbS \\
 e_5 &=& \mbH \\
 e_6 &=& \mbX
 \end{array}$
 & $q^{1 / 3} + y$
 & $\arraycolsep=1.4pt\begin{array}{rcl}
 \mbX &\leftrightarrow& -y \partial_x + p^2 \partial_p + 3 p q \partial_q - \tfrac{1}{2} y^2 \partial_z \\
 \mbY &\leftrightarrow& -(x \partial_y + \partial_p + \tfrac{1}{2} x^2 \partial_z) \\
 \mbH &\leftrightarrow& -x \partial_x + y \partial_y + 2 p \partial_p + 3 q \partial_q \qquad\qquad \\
 \mbU &\leftrightarrow& \partial_z \\
 \mbS &\leftrightarrow& \partial_x \\
 \mbT &\leftrightarrow& \partial_y + x \partial_z \\
 \end{array}$ \\
 \midrule
 $\mathbf{D.6}_{\lambda}$
 & $\begin{array}{c}
 \mfsl_2(\bbC) \oplus \mfsl_2(\bbC) \\
 \begin{array}{c|cccccc}
 & \mbX & \mbY & \mbH & \mbX' & \mbY' & \mbH' \\
	 \hline
 \mbX & \cdot & \mbH & -2 \mbX & \cdot & \cdot & \cdot \\
 \mbY & & \cdot & 2 \mbY & \cdot & \cdot & \cdot \\
 \mbH & & & \cdot & \cdot & \cdot & \cdot \\
 \mbX' & & & & \cdot & \mbH' & -2 \mbX' \\
 \mbY' & & & & & \cdot & 2 \mbY' \\
 \mbH' & & & & & & \cdot \\
 \end{array}
 \end{array}$
 & $\begin{array}{rl}
 \mfk\colon
 &\!\!\!\!\!\! \mbH + \mbH' \\
	 \mfd / \mfk\colon
 &\!\!\!\!\!\! \mbX - \lambda \mbX' , \\
 &\!\!\!\!\!\! \mbY - \mbY'
 \end{array}$
 & $\arraycolsep=1.4pt\begin{array}{rcl}
 e_1 &=& \mbX - \lambda^2 \mbX' \\
 e_2 &=& \mbY - \lambda \mbY' \\
 e_3 &=& \mbH + \lambda \mbH' \\
 e_4 &=& \mbX - \lambda \mbX' \\
 e_5 &=& \mbY - \mbY' \\
 e_6 &=& \mbH + \mbH'
 \end{array}$
 & $\begin{array}{c}
 y^{3 \alpha - 2} q^{\alpha} \\
 \lambda = (2 \alpha - 1)^2
 \end{array}$
 & $\arraycolsep=1.4pt\begin{array}{rcl}
 \mbX\phantom{'} &\leftrightarrow& \partial_x \\
 \mbY\phantom{'} &\leftrightarrow& -x^2 \partial_x - x y \partial_y + (x p - y) \partial_p + 3 x q \partial_q \\
 \mbH\phantom{'} &\leftrightarrow& -2 x \partial_x - y \partial_y + p \partial_p + 3 q \partial_q \\
 \mbX' &\leftrightarrow& \partial_z \\
 \mbY' &\leftrightarrow& \tfrac{1}{2 \alpha - 1} [\alpha y^{3 \alpha - 1} q^{\alpha - 1} \partial_x \\
 & & \quad + (\alpha y^{3 \alpha - 1} p q^{\alpha - 1} - y z) \partial_y \\
			 & & \quad + [(\alpha - 1) y^{3 \alpha - 1} q^{\alpha} - pz] \partial_p \\
 & & \quad + (-3 \alpha y^{3 \alpha - 2} p q^{\alpha} + q z) \partial_q \\
			 & & \quad + [\alpha (\alpha - 1) y^{6 \alpha - 3} q^{2 \alpha - 1} - (2 \alpha - 1)^2 z^2] \partial_z] \\
 \mbH' &\leftrightarrow& -\tfrac{1}{2 \alpha - 1}[y \partial_y + p \partial_p + q \partial_q + 2 (2 \alpha - 1) z \partial_z] \\
 \end{array}$ \\
 \midrule
 $\mathbf{D.6}_{\infty}$
 & $\begin{array}{c}
 \mfsl_2(\bbC) \oplus (\mfso_2(\bbC) \rightthreetimes \bbC^2) \\
 \begin{array}{c|cccccc}
 & \mbX & \mbY & \mbH & \mbZ & \mbV_1 & \mbV_2 \\
	 \hline
 \mbX & \cdot & \mbH & -2 \mbX & \cdot & \cdot & \cdot \\
 \mbY & & \cdot & 2 \mbY & \cdot & \cdot & \cdot \\
 \mbH & & & \cdot & \cdot & \cdot & \cdot \\
 \mbZ & & & & \cdot & \mbV_1 & -\mbV_2 \\
 \mbV_1 & & & & & \cdot & \cdot \\
 \mbV_2 & & & & & & \cdot \\
 \end{array}
 \end{array}$
 & $\begin{array}{rl}
 \mfk\colon
 &\!\!\!\!\!\! \mbH + 2 \mbZ \\
	 \mfd / \mfk\colon
 &\!\!\!\!\!\! \mbX + \mbV_1 , \\
 &\!\!\!\!\!\! \mbY + \mbV_2
 \end{array}$
 & $\arraycolsep=1.4pt\begin{array}{rcl}
 e_1 &=& \mbX \\
 e_2 &=& \mbY \\
 e_3 &=& \mbH \\
 e_4 &=& \mbX + \mbV_1 \\
 e_5 &=& \mbY + \mbV_2 \\
 e_6 &=& \mbH + 2 \mbZ
 \end{array}$
 & $p^{-2} q^2$
 & $\arraycolsep=1.4pt\begin{array}{rcl}
 \mbX &\leftrightarrow& -[\tfrac{1}{2} y^2 \partial_y + y p \partial_p + (y q + p^2) \partial_q + 2 p \partial_z] \\
 \mbY &\leftrightarrow& 2 \partial_y \\
 \mbH &\leftrightarrow& 2 ( y \partial_y + p \partial_p + q \partial_q ) \\
 \mbZ &\leftrightarrow& x \partial_x - p \partial_p - 2 q \partial_q - z \partial_z \\
 \mbV_1 &\leftrightarrow& 2 \partial_z \\
 \mbV_2 &\leftrightarrow& 2 \partial_x
 \end{array}$ \\
 \midrule
 $\mathbf{D.6}_{\ast}$
 & $\begin{array}{c}
 \mfso_3(\bbC) \rightthreetimes \bbC^3 \\
 \begin{array}{c|cccccc}
 & \mbX & \mbY & \mbH & \mbW_1 & \mbW_2 & \mbW_3 \\
	 \hline
 \mbX & \cdot & \mbH & -2 \mbX & \cdot & -\mbW_3 & 2 \mbW_1 \\
 \mbY & & \cdot & 2 \mbY & \mbW_3 & \cdot & -2 \mbW_2 \\
 \mbH & & & \cdot & 2 \mbW_1 & -2 \mbW_2 & \cdot \\
 \mbW_1 & & & & \cdot & \cdot & \cdot \\
 \mbW_2 & & & & & \cdot & \cdot \\
 \mbW_3 & & & & & & \cdot \\
 \end{array}
 \end{array}$
 & $\begin{array}{rl}
 \mfk\colon
 &\!\!\!\!\!\! \mbX + \mbY \\
	 \mfd / \mfk\colon
 &\!\!\!\!\!\! \mbH - 2 \mbW_3 ,\\
 &\!\!\!\!\!\! \mbX - \mbY \\
 &\!\!\!\!\!\! \quad + 2 \mbW_1 - 2 \mbW_2
 \end{array}$
 & $\arraycolsep=1.4pt\begin{array}{rcl}
 e_1 &=& -\mbW_1 + \mbW_2 \\
 e_2 &=& \mbW_3 \\
 e_3 &=& \mbW_1 + \mbW_2 \\
 e_4 &=& \mbH - 2 \mbW_3 \\
 e_5 &=& \mbX - \mbY \\
 & & \quad + 2 \mbW_1 - 2 \mbW_2 \\
 e_6 &=& \mbX + \mbY
 \end{array}$
 & $p^2 + q \log q$
 & $\arraycolsep=1.4pt\begin{array}{rcl}
 \mbX &\leftrightarrow& (2 y - \log q) \partial_x + [-p (\log q + 1) + z] \partial_y \\
 & & \quad - (p^2 + q) \partial_p - 4 p q \partial_q \\
 & & \quad - [(p^2 + q) \log q + 3 p^2 - q] \partial_z \\
 \mbY &\leftrightarrow& x \partial_y + \partial_p + (2 y - 1) \log z \\
 \mbH &\leftrightarrow& 2 [-x \partial_x + \partial_y + p \partial_p + 2 q \partial_q + (2 p + z) \partial_z] \\
 \mbW_1 &\leftrightarrow& - \partial_x\\
 \mbW_2 &\leftrightarrow& \partial_z\\
 \mbW_3 &\leftrightarrow& \partial_y
 \end{array}$ \\
 \bottomrule
\end{tabular}
\begin{minipage}{200mm}
\tiny
\vspace{0.25cm}
In the models $\mathbf{D.6}_{\lambda}$, $\lambda \in \bbC - \{0, \frac{1}{9}, 1, 9\}$ , and the models $\mathbf{D.6}_{\lambda}$ and $\mathbf{D.6}_{\lambda'}$ are equivalent iff $\lambda' = \lambda$ or $\lambda' = 1 / \lambda$.
\end{minipage}
\end{table}

\vspace*{\fill}

\end{landscape}

\subsection[Multiply transitive homogeneous real $(2, 3, 5)$ distributions]{Multiply transitive homogeneous real $\boldsymbol{(2, 3, 5)}$ distributions}

\begin{table}[h!]
\caption{The real form $\mbO^{\bbR}$ of model $\mbO$.}\label{table:real-O}\vspace{1mm}

\tiny
\centering
\begin{tabular}{lcl}
\toprule
 real form
 & symmetry algebra
 & \multicolumn{1}{c}{algebraic model} \\
\midrule
 $\mbO^{\bbR}$
 & $\mfg_2$
 & $\begin{array}{rl}
 \mfk\colon
 &\!\!\!\!\!\! \mbA_{11}, \mbA_{12}, \mbA_{21}, \mbA_{22}, \mbZ_1, \mbZ_2, \mbs, \mbW_1, \mbW_2 \\
	 \mfd / \mfk\colon
 &\!\!\!\!\!\! \mbX_1, \mbX_2
 \end{array}$ \\
\bottomrule
\end{tabular}

\vspace{0.8cm}

\caption{Real forms of the models $\mathbf{N.7}_{\Lambda}$.}\label{table:real-N.7}\vspace{1mm}

\tiny
\centering
\begin{tabular}{cclcl}
\toprule
 $\Lambda$
 & real form
 & \multicolumn{1}{c}{$r, s, a, b$}
 & $\mbF$-action
 & \multicolumn{1}{c}{algebraic model} \\
\midrule
 \multirow{5.5}{*}{$(-\infty, 0) \cup (1, \infty)$}
 & $\mathbf{N.7}_{\Lambda}^{\nwarrow}$
	& $\left\{
		\begin{aligned}
			r &= -a^2 - b^2 \\
			s &= a^2 b^2
		\end{aligned} \right.$
 & $\begin{pmatrix}
 \cdot & \cdot & -a & \cdot \\
 \cdot & \cdot & \cdot & b \\
 a & \cdot & \cdot & \cdot \\
 \cdot & -b & \cdot & \cdot
 \end{pmatrix}$
 & $\begin{array}{rl}
 \mfk\colon
 &\!\!\!\!\!\! \mbE, a^{-1 / 2} \mbS_1 + b^{-1 / 2} \mbS_2 \\
	 \mfd / \mfk\colon
 &\!\!\!\!\!\! \mbF, -a^{ 1 / 2} \mbT_1 + b^{ 1 / 2} \mbT_2
 \end{array}$ \\
 \cmidrule{2-5}
 & $\mathbf{N.7}_{\Lambda}^{\nearrow}$
	& $\left\{
		\begin{aligned}
			r &= a^2 + b^2 \\
			s &= a^2 b^2
		\end{aligned} \right.$
 & $\begin{pmatrix}
 a & \cdot & \cdot & \cdot \\
 \cdot & b & \cdot & \cdot \\
 \cdot & \cdot & -a & \cdot \\
 \cdot & \cdot & \cdot & -b
 \end{pmatrix}$
 & $\begin{array}{rl}
 \mfk\colon
 &\!\!\!\!\!\! \mbE, \tfrac{1}{a} \mbS_1 + \tfrac{1}{b} \mbS_2 - \mbT_1 + \mbT_2 \\
	 \mfd / \mfk\colon
 &\!\!\!\!\!\! \mbF, \mbS_1 + \mbS_2 + a \mbT_1 - b \mbT_2
 \end{array}$ \\
 \cmidrule{1-5}
 \multirow{5.5}{*}{$0$}
 & $\mathbf{N.7}_0^{\leftarrow}$
	& $\left\{
		\begin{aligned}
			r &= -1 {}^{\phantom{2}} \\
			s &= 0 {}^{\phantom{2}}
		\end{aligned} \right.$
 & $\begin{pmatrix}
 \cdot & \cdot & -1 & \cdot \\
 \cdot & \cdot & \cdot & 1 \\
 1 & \cdot & \cdot & \cdot \\
 \cdot & \cdot & \cdot & \cdot
 \end{pmatrix}$
 & $\begin{array}{rl}
 \mfk\colon
 &\!\!\!\!\!\! \mbT_1 + \mbT_2 \\
	 \mfd / \mfk\colon
 &\!\!\!\!\!\! \mbS_1 - \mbS_2
 \end{array}$ \\
 \cmidrule{2-5}
 & $\mathbf{N.7}_0^{\rightarrow}$
	& $\left\{
		\begin{aligned}
			r &= 1 {}^{\phantom{2}} \\
			s &= 0 {}^{\phantom{2}}
		\end{aligned} \right.$
 & $\begin{pmatrix}
 1 & \cdot & \cdot & \cdot \\
 \cdot & \cdot & \cdot & 1 \\
 \cdot & \cdot & -1 & \cdot \\
 \cdot & \cdot & \cdot & \cdot
 \end{pmatrix}$
 & $\begin{array}{rl}
 \mfk\colon
 &\!\!\!\!\!\! \mbE, -2 \mbS_1 + \mbT_1 + 2 \mbT_2 \\
	 \mfd / \mfk\colon
 &\!\!\!\!\!\! \mbF, - \mbS_2 + \mbT_1 + \mbT_2
 \end{array}$ \\
 \cmidrule{1-5}
 \multirow{5.5}{*}{$(0, \frac{16}{25})$}
 & $\mathbf{N.7}_\Lambda^\swarrow$
	& $\left\{
		\begin{aligned}
			r &= a^2 - b^2 < 0 \\
			s &= -a^2 b^2
		\end{aligned} \right.$
 & $\begin{pmatrix}
 a & \cdot & \cdot & \cdot \\
 \cdot & \cdot & \cdot & -b \\
 \cdot & \cdot & -a & \cdot \\
 \cdot & b & \cdot & \cdot
 \end{pmatrix}$
 & $\begin{array}{rl}
 \mfk\colon
 &\!\!\!\!\!\! \mbE, a \mbS_1 + 2 a \mbS_2 + 2 b \mbT_1 \\
	 \mfd / \mfk\colon
 &\!\!\!\!\!\! \mbF, a \mbS_1 + a \mbS_2 + b \mbT_2
 \end{array}$ \\
 \cmidrule{2-5}
 & $\mathbf{N.7}_\Lambda^\searrow$
	& $\left\{
		\begin{aligned}
			r &= a^2 - b^2 > 0 \\
			s &= -a^2 b^2
		\end{aligned} \right.$
 & $\begin{pmatrix}
 a & \cdot & \cdot & \cdot \\
 \cdot & \cdot & \cdot & -b \\
 \cdot & \cdot & -a & \cdot \\
 \cdot & b & \cdot & \cdot
 \end{pmatrix}$
 & $\begin{array}{rl}
 \mfk\colon
 &\!\!\!\!\!\! \mbE, a \mbS_1 + 2 a \mbS_2 + 2 b \mbT_1 \\
	 \mfd / \mfk\colon
 &\!\!\!\!\!\! \mbF, a \mbS_1 + a \mbS_2 + b \mbT_2
 \end{array}$ \\
 \cmidrule{1-5}
 \multirow{5.5}{*}{$\frac{16}{25}$}
 & $\mathbf{N.7}_{16 / 25}^{\downarrow}$
	& $\left\{
		\begin{aligned}
			r &= 0 {}^{\phantom{2}} \\
			s &= -1 {}^{\phantom{2}}
		\end{aligned} \right.$
 & $\begin{pmatrix}
 1 & \cdot & \cdot & \cdot \\
 \cdot & \cdot & \cdot & -1 \\
 \cdot & \cdot & -1 & \cdot \\
 \cdot & 1 & \cdot & \cdot
 \end{pmatrix}$
 & $\begin{array}{rl}
 \mfk\colon
 &\!\!\!\!\!\! \mbE, \mbS_1 + 2 \mbS_2 + 2 \mbT_1 \\
	 \mfd / \mfk\colon
 &\!\!\!\!\!\! \mbF, \mbS_1 + \mbS_2 + \mbT_2
 \end{array}$ \\
 \cmidrule{2-5}
 & $\mathbf{N.7}_{16 / 25}^{\uparrow}$
	& $\left\{
		\begin{aligned}
			r &= 0 {}^{\phantom{2}} \\
			s &= 1 {}^{\phantom{2}}
		\end{aligned} \right.$
 & $\frac{1}{\sqrt{2}}\begin{pmatrix}
 1 & \cdot & \cdot & -1 \\
 \cdot & -1 & -1 & \cdot \\
 \cdot & 1 & -1 & \cdot \\
 1 & \cdot & \cdot & 1
 \end{pmatrix}$
 & $\begin{array}{rl}
 \mfk\colon
 &\!\!\!\!\!\! \mbE, \mbS_2 - \mbT_1 - 2 \mbT_2 \\
	 \mfd / \mfk\colon
 &\!\!\!\!\!\! \mbF, \mbS_1 + \mbT_1 - \mbT_2
 \end{array}$ \\
 \cmidrule{1-5}
 \multirow{5.5}{*}{$(\frac{16}{25}, 1)$}
 & $\mathbf{N.7}_\Lambda^\nwarrow$
	& $\left\{
		\begin{aligned}
			r &= 2 \big(a^2 - b^2\big) < 0 \\
			s &= \big(a^2 + b^2\big)^2
		\end{aligned} \right.$
 & $\begin{pmatrix}
 a & \cdot & \cdot & -b \\
 \cdot & -a & -b & \cdot \\
 \cdot & b & -a & \cdot \\
 b & \cdot & \cdot & a
 \end{pmatrix}$
 & $\begin{array}{rl}
 \mfk\colon
 &\!\!\!\!\!\! \mbE, ab \mbS_2 - a^2 \mbT_1 - \big(a^2 + b^2\big) \mbT_2 \\
	 \mfd / \mfk\colon
 &\!\!\!\!\!\! \mbF, b \mbS_1 + a \mbT_1 - a \mbT_2
 \end{array}$ \\
 \cmidrule{2-5}
 & $\mathbf{N.7}_\Lambda^\nearrow$
	& $\left\{
		\begin{aligned}
			r &= 2 \big(a^2 - b^2\big) > 0 \\
			s &= \big(a^2 + b^2\big)^2
		\end{aligned} \right.$
 & $\begin{pmatrix}
 a & \cdot & \cdot & -b \\
 \cdot & -a & -b & \cdot \\
 \cdot & b & -a & \cdot \\
 b & \cdot & \cdot & a
 \end{pmatrix}$
 & $\begin{array}{rl}
 \mfk\colon
 &\!\!\!\!\!\! \mbE, ab \mbS_2 - a^2 \mbT_1 - \big(a^2 + b^2\big) \mbT_2 \\
	 \mfd / \mfk\colon
 &\!\!\!\!\!\! \mbF, b \mbS_1 + a \mbT_1 - a \mbT_2
 \end{array}$ \\
 \cmidrule{1-5}
 \multirow{5.5}{*}{$1$}
 & $\mathbf{N.7}_1^{\nwarrow}$
	& $\left\{
		\begin{aligned}
			r &= -2 {}^{\phantom{2}} \\
			s &= 1 {}^{\phantom{2}}
		\end{aligned} \right.$
 & $\begin{pmatrix}
 \cdot & -1 & 1 & \cdot \\
 1 & \cdot & \cdot & 1 \\
 \cdot & \cdot & \cdot & -1 \\
 \cdot & \cdot & 1 & \cdot
 \end{pmatrix}$
 & $\begin{array}{rl}
 \mfk\colon
 &\!\!\!\!\!\! \mbE, \mbS_2 + 2 \mbT_1 \\
	 \mfd / \mfk\colon
 &\!\!\!\!\!\! \mbF, \mbS_1 + 2 \mbT_2
 \end{array}$ \\
 \cmidrule{2-5}
 & $\mathbf{N.7}_1^{\nearrow}$
	& $\left\{
		\begin{aligned}
			r &= 2 {}^{\phantom{2}} \\
			s &= 1 {}^{\phantom{2}}
		\end{aligned} \right.$
 & $\begin{pmatrix}
 1 & \cdot & \cdot & 1 \\
 \cdot & -1 & 1 & \cdot \\
 \cdot & \cdot & -1 & \cdot \\
 \cdot & \cdot & \cdot & 1
 \end{pmatrix}$
 & $\begin{array}{rl}
 \mfk\colon
 &\!\!\!\!\!\! \mbE, \mbS_2 + \mbT_1 + \mbT_2 \\
	 \mfd / \mfk\colon
 &\!\!\!\!\!\! \mbF, \mbS_1 + \mbS_2 + 2 \mbT_2
 \end{array}$ \\
\bottomrule
\end{tabular}
\begin{minipage}{113mm}
\tiny
\vspace{0.25cm}
For each $\Lambda \in \bbR$, $(r, s)$ is a pair that satisfies $\Lambda = 64 s / \big(100 s - 9 r^2\big)$; see Section~\ref{subsection:N.7-real} for details.
\end{minipage}
\end{table}

\newpage

\begin{table}[th!]
\caption{Real forms of model $\mathbf{N.6}$.}\label{table:real-N.6}\vspace{1mm}

\tiny
\centering
\begin{tabular}{lccl}
\toprule
 real form
 & anti-involution
 & symmetry algebra structure
 & \multicolumn{1}{c}{algebraic model} \\
\midrule
 $\mathbf{N.6}^+$
 & $\diag(1,1,1,1,1,1)$
 &
 $\begin{array}{c}
 \mfsl_2 \rightthreetimes \mfm \\
 \begin{array}{c|cccccc}
 & \mbX & \mbY & \mbH & \mbU & \mbS & \mbT \\
	 \hline
 \mbX & \cdot & \mbH & -2 \mbX & \cdot & \cdot & \mbS \\
 \mbY & & \cdot & 2 \mbY & \cdot & \mbT & \cdot \\
 \mbH & & & \cdot & \cdot & \mbS & -\mbT \\
 \mbU & & & & \cdot & \cdot & \cdot \\
 \mbS & & & & & \cdot & \mbU \\
 \mbT & & & & & & \cdot \\
 \end{array}
 \end{array}$
 & $\begin{array}{rl}
 \mfk\colon
 &\!\!\!\!\!\! \mbX \\
	 \mfd / \mfk\colon
 &\!\!\!\!\!\! \mbY - \mbU - \mbS, \mbH
 \end{array}$ \\
\cmidrule{1-4}
 $\mathbf{N.6}^-$
 & $\diag(1, -1, -1, -1, 1, -1)$
 &
 $\begin{array}{c}
 \mfsl_2 \rightthreetimes \mfm \\
 \begin{array}{c|cccccc}
 & \mbX & \mbY & \mbH & \mbU & \mbS & \mbT \\
	 \hline
 \mbX & \cdot & \mbH & -2 \mbX & \cdot & \cdot & \mbS \\
 \mbY & & \cdot & 2 \mbY & \cdot & \mbT & \cdot \\
 \mbH & & & \cdot & \cdot & \mbS & -\mbT \\
 \mbU & & & & \cdot & \cdot & \cdot \\
 \mbS & & & & & \cdot & \mbU \\
 \mbT & & & & & & \cdot \\
 \end{array}
 \end{array}$
 & $\begin{array}{rl}
 \mfk\colon
 &\!\!\!\!\!\! \mbX \\
	 \mfd / \mfk\colon
 &\!\!\!\!\!\! \mbY + \mbU + \mbS, \mbH
 \end{array}$ \\
\bottomrule
\end{tabular}

\vspace{0.8cm}

\caption{Real forms of the models $\mathbf{D.6}_{\lambda}$.}\label{table:real-D.6-lambda}\vspace{1mm}

\tiny
\centering
\setlength\tabcolsep{0.9pt}
\begin{tabular}{lccccl}
 \toprule
 real form
 & $\lambda$
 & anti-involution
 & symmetry algebra structure
 & \multicolumn{1}{c}{algebraic model} \\
 \midrule
 $\mathbf{D.6}_{\lambda}^{2-}$
 & $\begin{array}{c}
 \bbR - \{0, \frac{1}{9}, 1, 9\} \\
 \lambda \sim 1 / \lambda
 \end{array}$
 & $\diag(1,1,1,1,1,1)$
 & $\begin{array}{c}
 \mfsl_2 \oplus \mfsl_2 \\
 \begin{array}{c|cccccc}
 & \mbX & \mbY & \mbH & \mbX' & \mbY' & \mbH' \\
	 \hline
 \mbX & \cdot & \mbH & -2 \mbX & \cdot & \cdot & \cdot \\
 \mbY & & \cdot & 2 \mbY & \cdot & \cdot & \cdot \\
 \mbH & & & \cdot & \cdot & \cdot & \cdot \\
 \mbX' & & & & \cdot & \mbH' & -2 \mbX' \\
 \mbY' & & & & & \cdot & 2 \mbY' \\
 \mbH' & & & & & & \cdot \\
 \end{array}
 \end{array}$
 & $\begin{array}{rl}
 \mfk\colon
 &\!\!\!\!\!\! \mbH + \mbH' \\
	 \mfd / \mfk\colon
 &\!\!\!\!\!\! \mbX - \lambda \mbX', \\
 &\!\!\!\!\!\! \mbY - \mbY'
 \end{array}$ \\
 \cmidrule{1-5}
 $\mathbf{D.6}_{\lambda}^{2+}$
 & $\begin{array}{c}
 \bbR^+ - \{\frac{1}{9}, 1, 9\} \\
 \lambda \sim 1 / \lambda
 \end{array}$
 & $\begin{pmatrix}
 \cdot & 1 & \cdot & \cdot & \cdot & \cdot \\
 1 & \cdot & \cdot & \cdot & \cdot & \cdot \\
 \cdot & \cdot & -1 & \cdot & \cdot & \cdot \\
 \cdot & \cdot & \cdot & \cdot & 1 & \cdot \\
 \cdot & \cdot & \cdot & 1 & \cdot & \cdot \\
 \cdot & \cdot & \cdot & \cdot & \cdot & -1 \\
 \end{pmatrix}$
 & $\begin{array}{c}
 \mfsl_2 \oplus \mfsl_2 \\
 \begin{array}{c|cccccc}
 & \mbX & \mbY & \mbH & \mbX' & \mbY' & \mbH' \\
	 \hline
 \mbX & \cdot & \mbH & -2 \mbX & \cdot & \cdot & \cdot \\
 \mbY & & \cdot & 2 \mbY & \cdot & \cdot & \cdot \\
 \mbH & & & \cdot & \cdot & \cdot & \cdot \\
 \mbX' & & & & \cdot & \mbH' & -2 \mbX' \\
 \mbY' & & & & & \cdot & 2 \mbY' \\
 \mbH' & & & & & & \cdot \\
 \end{array}
 \end{array}$
 & $\begin{array}{rl}
 \mfk\colon
 &\!\!\!\!\!\! \wt\mbH + \wt\mbH' \\
	 \mfd / \mfk\colon
 &\!\!\!\!\!\! \wt\mbX + \sqrt{\lambda} \wt\mbX', \\
 &\!\!\!\!\!\! \wt\mbY + \sqrt{\lambda} \wt\mbY'
 \end{array}$ \\
 \cmidrule{1-5}
 $\mathbf{D.6}_{\lambda}^4$
 & $\begin{array}{c}
 (-\infty, 0)
 \end{array}$
 & $\begin{pmatrix}
 \cdot & 1 & \cdot & \cdot & \cdot & \cdot \\
 1 & \cdot & \cdot & \cdot & \cdot & \cdot \\
 \cdot & \cdot & -1 & \cdot & \cdot & \cdot \\
 \cdot & \cdot & \cdot & \cdot & 1 & \cdot \\
 \cdot & \cdot & \cdot & 1 & \cdot & \cdot \\
 \cdot & \cdot & \cdot & \cdot & \cdot & -1 \\
 \end{pmatrix}$
 & $\begin{array}{c}
 \mfsl_2\oplus \mfso_3 \\
 \begin{array}{c|cccccc}
 & \mbX & \mbY & \mbH & \mbA & \mbB & \mbC \\
	 \hline
 \mbX & \cdot & \mbH & -2 \mbX & \cdot & \cdot & \cdot \\
 \mbY & & \cdot & 2 \mbY & \cdot & \cdot & \cdot \\
 \mbH & & & \cdot & \cdot & \cdot & \cdot \\
 \mbA & & & & \cdot & \mbC & -\mbB \\
 \mbB & & & & & \cdot & \mbA \\
 \mbC & & & & & & \cdot \\
 \end{array}
 \end{array}$
 & $\begin{array}{rl}
 \mfk\colon
 &\!\!\!\!\!\! \wt\mbH + \mbC \\
	 \mfd / \mfk\colon
 &\!\!\!\!\!\! \wt\mbX + \sqrt{-\lambda} \mbA, \\
 &\!\!\!\!\!\! \wt\mbY + \sqrt{-\lambda} \mbB
 \end{array}$ \\
 \cmidrule{1-5}
 $\mathbf{D.6}_{\lambda}^6$
 & $\begin{array}{c}
 \bbR^+ - \{\frac{1}{9}, 1, 9\} \\
 \lambda \sim 1 / \lambda
 \end{array}$
 & $\begin{pmatrix}
 \cdot & -1 & \cdot & \cdot & \cdot & \cdot \\
 -1 & \cdot & \cdot & \cdot & \cdot & \cdot \\
 \cdot & \cdot & -1 & \cdot & \cdot & \cdot \\
 \cdot & \cdot & \cdot & \cdot & -1 & \cdot \\
 \cdot & \cdot & \cdot & -1 & \cdot & \cdot \\
 \cdot & \cdot & \cdot & \cdot & \cdot & -1 \\
 \end{pmatrix}$
 & $\begin{array}{c}
 \mfso_3 \oplus \mfso_3 \\
 \begin{array}{c|cccccc}
 & \mbA & \mbB & \mbC & \mbA' & \mbB' & \mbC' \\
	 \hline
 \mbA & \cdot & \mbC & -\mbB & \cdot & \cdot & \cdot \\
 \mbB & & \cdot & \mbA & \cdot & \cdot & \cdot \\
 \mbC & & & \cdot & \cdot & \cdot & \cdot \\
 \mbA' & & & & \cdot & \mbC' & -\mbB' \\
 \mbB' & & & & & \cdot & \mbA' \\
 \mbC' & & & & & & \cdot \\
 \end{array}
 \end{array}$
 & $\begin{array}{rl}
 \mfk\colon
 &\!\!\!\!\!\! \mbC + \mbC' \\
	 \mfd / \mfk\colon
 &\!\!\!\!\!\! \mbA + \sqrt{-\lambda} \mbA', \\
 &\!\!\!\!\!\! \mbB + \sqrt{-\lambda} \mbB'
 \end{array}$ \\
 \cmidrule{1-5}
 $\mathbf{D.6}_{-1}^{3-}$
 & $-1$
 & $\diag(-1, 1, -1, 1, -1, 1)$
 &
 $\begin{array}{c}
 \mfso_{1, 3} \\
 \begin{array}{c|cccccc}
		 & \mbA & \mbB & \mbC & \mbD_{\mbA} & \mbD_{\mbB} & \mbD_{\mbC} \\
 \hline
	 \mbA & \cdot & \mbC & -\mbB & \cdot & \mbD_{\mbC} & -\mbD_{\mbB} \\
		 \mbB & & \cdot & \mbA & -\mbD_{\mbC} & \cdot & \mbD_{\mbA} \\
		 \mbC & & & \cdot & \mbD_{\mbB} & -\mbD_{\mbA} & \cdot \\
		 \mbD_{\mbA} & & & & \cdot & \mbC & -\mbB \\
		 \mbD_{\mbB} & & & & & \cdot & \mbA \\
		 \mbD_{\mbC} & & & & & & \cdot \\
	 \end{array}
	 \end{array}$
 & $\begin{array}{rl}
 \mfk\colon
 &\!\!\!\!\!\! \mbD_\mbC \\
	 \mfd / \mfk\colon
 &\!\!\!\!\!\! \mbA + \mbD_\mbB , \\
 &\!\!\!\!\!\! \mbB + \mbD_\mbA
 \end{array}$ \\
 \cmidrule{1-5}
 $\mathbf{D.6}_{-1}^{3+}$
 & $-1$
 & $\begin{pmatrix}
 \cdot & -i & \cdot & \cdot & \cdot & \cdot \\
 -i & \cdot & \cdot & \cdot & \cdot & \cdot \\
 \cdot & \cdot & 1 & \cdot & \cdot & \cdot \\
 \cdot & \cdot & \cdot & \cdot & i & \cdot \\
 \cdot & \cdot & \cdot & i & \cdot & \cdot \\
 \cdot & \cdot & \cdot & \cdot & \cdot & -1 \\
 \end{pmatrix}$
 &
 $\begin{array}{c}
 \mfso_{1, 3} \\
 \begin{array}{c|cccccc}
		 & \mbA & \mbB & \mbC & \mbD_{\mbA} & \mbD_{\mbB} & \mbD_{\mbC} \\
 \hline
	 \mbA & \cdot & \mbC & -\mbB & \cdot & \mbD_{\mbC} & -\mbD_{\mbB} \\
		 \mbB & & \cdot & \mbA & -\mbD_{\mbC} & \cdot & \mbD_{\mbA} \\
		 \mbC & & & \cdot & \mbD_{\mbB} & -\mbD_{\mbA} & \cdot \\
		 \mbD_{\mbA} & & & & \cdot & \mbC & -\mbB \\
		 \mbD_{\mbB} & & & & & \cdot & \mbA \\
		 \mbD_{\mbC} & & & & & & \cdot \\
	 \end{array}
	 \end{array}$
 & $\begin{array}{rl}
 \mfk\colon
 &\!\!\!\!\!\! \mbC \\
	 \mfd / \mfk\colon
 &\!\!\!\!\!\! \mbA + \mbD_\mbA , \\
 &\!\!\!\!\!\! \mbB + \mbD_\mbB
 \end{array}$ \\
 \bottomrule
\end{tabular}
\begin{minipage}{113mm}
\tiny
\vspace{0.25cm}
In models $\mathbf{D.6}_{\lambda}^{2+}$ and $\mathbf{D.6}_{\lambda}^4$,
 $\wt\mbX := \tfrac{1}{\sqrt{2}} \mbX + \tfrac{1}{2} \mbH$ ,
 $\wt\mbY := -\tfrac{1}{\sqrt{2}} \mbY + \tfrac{1}{2} \mbH$ ,
 $\wt\mbH := \tfrac{1}{\sqrt{2}} (\mbX - \mbY) + \tfrac{1}{2} \mbH$ .
\end{minipage}
\end{table}

\newpage

\begin{table}[th!]
\caption{Real forms of model $\mathbf{D.6}_\infty$.}\label{table:real-D.6-infty}\vspace{1mm}

\tiny
\centering
\begin{tabular}{lccl}
\toprule
 real form
 & anti-involution
 & symmetry algebra structure
 & \multicolumn{1}{c}{algebraic model} \\
\midrule
 $\mathbf{D.6}_\infty^1$
 & $\diag(1,1,1,1,1,1)$
 &
 $\begin{array}{c}
 \mfsl_2 \oplus \big(\mfso_{1, 1} \rightthreetimes \bbR^{1, 1}\big) \\
 \begin{array}{c|cccccc}
 & \mbX & \mbY & \mbH & \mbZ & \mbV_1 & \mbV_2 \\
	 \hline
 \mbX & \cdot & \mbH & -2 \mbX & \cdot & \cdot & \cdot \\
 \mbY & & \cdot & 2 \mbY & \cdot & \cdot & \cdot \\
 \mbH & & & \cdot & \cdot & \cdot & \cdot \\
 \mbZ & & & & \cdot & \mbV_1 & -\mbV_2 \\
 \mbV_1 & & & & & \cdot & \cdot \\
 \mbV_2 & & & & & & \cdot \\
 \end{array}
 \end{array}$
 & $\begin{array}{rl}
 \mfk\colon
 &\!\!\!\!\!\! \mbX - \mbY + 2 \mbZ \\
	 \mfd / \mfk\colon
 &\!\!\!\!\!\! \mbX + \mbY - \mbH, \\
 &\!\!\!\!\!\!
 \end{array}$ \\
\cmidrule{1-4}
 $\mathbf{D.6}_\infty^2$
 & $\begin{pmatrix}
 \cdot & 1 & \cdot & \cdot & \cdot & \cdot \\
 1 & \cdot & \cdot & \cdot & \cdot & \cdot \\
 \cdot & \cdot & -1 & \cdot & \cdot & \cdot \\
 \cdot & \cdot & \cdot & \cdot & 1 & \cdot \\
 \cdot & \cdot & \cdot & 1 & \cdot & \cdot \\
 \cdot & \cdot & \cdot & \cdot & \cdot & -1 \\
 \end{pmatrix}$
 &
 $\begin{array}{c}
 \mfsl_2 \oplus \big(\mfso_2 \rightthreetimes \bbR^2\big) \\
 \begin{array}{c|cccccc}
 & \mbX & \mbY & \mbH & \hat\mbZ & \mbV_1 & \mbV_2 \\
	 \hline
 \mbX & \cdot & \mbH & -2 \mbX & \cdot & \cdot & \cdot \\
 \mbY & & \cdot & 2 \mbY & \cdot & \cdot & \cdot \\
 \mbH & & & \cdot & \cdot & \cdot & \cdot \\
 \hat\mbZ & & & & \cdot & \mbV_2 & -\mbV_1 \\
 \mbV_1 & & & & & \cdot & \cdot \\
 \mbV_2 & & & & & & \cdot \\
 \end{array}
 \end{array}$
 & $\begin{array}{rl}
 \mfk\colon
 &\!\!\!\!\!\! \mbX - \mbY + 2 \hat\mbZ \\
	 \mfd / \mfk\colon
 &\!\!\!\!\!\! \sqrt{2} \mbX + \sqrt{2} \mbY - \mbV_1 - \mbV_2 , \\
			&\!\!\!\!\!\! \sqrt{2} \mbH + \mbV_1 - \mbV_2
 \end{array}$ \\
\cmidrule{1-4}
 $\mathbf{D.6}_\infty^4$
 & $\begin{pmatrix}
 \cdot & -1 & \cdot & \cdot & \cdot & \cdot \\
 -1 & \cdot & \cdot & \cdot & \cdot & \cdot \\
 \cdot & \cdot & -1 & \cdot & \cdot & \cdot \\
 \cdot & \cdot & \cdot & \cdot & -1 & \cdot \\
 \cdot & \cdot & \cdot & -1 & \cdot & \cdot \\
 \cdot & \cdot & \cdot & \cdot & \cdot & -1 \\
 \end{pmatrix}$
 &
 $\begin{array}{c}
 \mfso_3 \oplus \big(\mfso_2 \rightthreetimes \bbR^2\big) \\
 \begin{array}{c|cccccc}
 & \mbA & \mbB & \mbC & \hat\mbZ & \mbV_1 & \mbV_2 \\
	 \hline
 \mbA & \cdot & \mbC & -\mbB & \cdot & \cdot & \cdot \\
 \mbB & & \cdot & \mbA & \cdot & \cdot & \cdot \\
 \mbC & & & \cdot & \cdot & \cdot & \cdot \\
 \hat\mbZ & & & & \cdot & \mbV_2 & -\mbV_1 \\
 \mbV_1 & & & & & \cdot & \cdot \\
 \mbV_2 & & & & & & \cdot \\
 \end{array}
 \end{array}$
 & $\begin{array}{rl}
 \mfk\colon
 &\!\!\!\!\!\! \mbC + \hat\mbZ \\
	 \mfd / \mfk\colon
 &\!\!\!\!\!\! \mbA + \mbV_1 , \\
			&\!\!\!\!\!\! \mbB + \mbV_2
 \end{array}$ \\
\bottomrule
\end{tabular}

\vspace{10mm}

\caption{Real forms of model $\mathbf{D.6}_\ast$.}\label{table:real-D.6-ast}\vspace{1mm}

\tiny
\centering
\begin{tabular}{lccl}
\toprule
 real form
 & anti-involution
 & symmetry algebra structure
 & \multicolumn{1}{c}{algebraic model} \\
\midrule
 $\mathbf{D.6}_\ast^{1-}$
 & $\diag(1,1,1,1,1,1)$
 & $\begin{array}{c}
 \mfsl_2 \rightthreetimes \bbR^{1, 2} \\
 \begin{array}{c|cccccc}
 & \mbX & \mbY & \mbH & \mbW_1 & \mbW_2 & \mbW_3 \\
	 \hline
 \mbX & \cdot & \mbH & -2 \mbX & \cdot & -\mbW_3 & 2 \mbW_1 \\
 \mbY & & \cdot & 2 \mbY & \mbW_3 & \cdot & -2 \mbW_2 \\
 \mbH & & & \cdot & 2 \mbW_1 & -2 \mbW_2 & \cdot \\
 \mbW_1 & & & & \cdot & \cdot & \cdot \\
 \mbW_2 & & & & & \cdot & \cdot \\
 \mbW_3 & & & & & & \cdot \\
 \end{array}
 \end{array}$
 & $\begin{array}{rl}
 \mfk\colon
 &\!\!\!\!\!\! \mbX + \mbY \\
	 \mfd / \mfk\colon
 &\!\!\!\!\!\! \mbH - 2 \mbW_3 , \\
 &\!\!\!\!\!\! \mbX - \mbY + 2 \mbW_1 - 2 \mbW_2
 \end{array}$ \\
\cmidrule{1-4}
 $\mathbf{D.6}_\ast^{1+}$
 & $\diag(-1, 1, -1, 1, -1, -1)$
 & $\begin{array}{c}
 \mfsl_2 \rightthreetimes \bbR^{1, 2} \\
 \begin{array}{c|cccccc}
 & \mbX & \mbY & \mbH & \mbW_1 & \mbW_2 & \mbW_3 \\
	 \hline
 \mbX & \cdot & \mbH & -2 \mbX & \cdot & -\mbW_3 & 2 \mbW_1 \\
 \mbY & & \cdot & 2 \mbY & \mbW_3 & \cdot & -2 \mbW_2 \\
 \mbH & & & \cdot & 2 \mbW_1 & -2 \mbW_2 & \cdot \\
 \mbW_1 & & & & \cdot & \cdot & \cdot \\
 \mbW_2 & & & & & \cdot & \cdot \\
 \mbW_3 & & & & & & \cdot \\
 \end{array}
 \end{array}$
 & $\begin{array}{rl}
 \mfk\colon
 &\!\!\!\!\!\! \mbX - \mbY \\
	 \mfd / \mfk\colon
 &\!\!\!\!\!\! \mbH + 2 \mbW_3 , \\
			&\!\!\!\!\!\! \mbX + \mbY - 2 \mbW_1 - 2 \mbW_2
 \end{array}$ \\
\cmidrule{1-4}
 $\mathbf{D.6}_\ast^3$
 & $\diag(1, -1, -1, -1, 1, -1)$
 &
 $\begin{array}{c}
 \mfso_3 \rightthreetimes \bbR^3 \\
 \begin{array}{c|cccccc}
 & \mbA & \mbB & \mbC & \mbW_\mbA & \mbW_\mbB & \mbW_\mbC \\
	 \hline
 \mbA & \cdot & \mbC & -\mbB & \cdot & \mbW_\mbC & -\mbW_\mbB \\
 \mbB & & \cdot & \mbA & -\mbW_\mbC & \cdot & \mbW_\mbA \\
 \mbC & & & \cdot & \mbW_\mbB & -\mbW_\mbA & \cdot \\
 \mbW_\mbA & & & & \cdot & \cdot & \cdot \\
 \mbW_\mbB & & & & & \cdot & \cdot \\
 \mbW_\mbC & & & & & & \cdot \\
 \end{array}
 \end{array}$
 & $\begin{array}{rl}
 \mfk\colon
 &\!\!\!\!\!\! \mbC \\
	 \mfd / \mfk\colon
 &\!\!\!\!\!\! \mbA - \mbW_\mbA , \\
			&\!\!\!\!\!\! \mbB - \mbW_\mbB
 \end{array}$ \\
\bottomrule
\end{tabular}
\end{table}

\subsection*{Acknowledgements}
It is a pleasure to thank Boris Doubrov, both for several helpful conversations and (again) for access to the unpublished notes \cite{DoubrovGovorovClassification}. It is likewise a pleasure to thank Dennis The for an invaluable exchange about the classification of real forms of a given complex geometric structure. The author is also grateful to the two anonymous referees for several helpful comments and two small corrections. Finally, the author gratefully acknowledges support from the Austrian Science Fund (FWF), via project P27072--N25, the Simons Foundation, via grant 346300 for IMPAN, and the Polish Government, via the matching 2015-2019 Polish MNiSW fund.

\newpage
\pdfbookmark[1]{References}{ref}
\LastPageEnding

\end{document}